\documentclass[aps,physrev,preprint,groupedaddress]{revtex4-2}

\usepackage{amsmath}
\usepackage{txfonts}
\usepackage{graphicx}
\usepackage{dcolumn}
\usepackage{bm}
\usepackage{epsfig}
\usepackage{booktabs}
\usepackage{subfigure}
\usepackage{graphics}
\usepackage{amssymb}
\usepackage{array}
\usepackage{color}
\usepackage{tabularx}
\usepackage{multirow}

\usepackage{CJK}
\usepackage{xcolor}
\usepackage{algorithm}
\usepackage{algpseudocode}

\begin{document}


\title{Equilibrium-distribution-function based mesoscopic finite-difference methods for partial differential equations: Modeling and Analysis}



\author{Baochang Shi }
\email[This work is dedicated to the memory of Prof. Nengchao Wang, my esteemed mentor, with deepest gratitude. Corresponding author] {(shibc@hust.edu.cn)}
\affiliation{School of Mathematics and Statistics, Huazhong
University of Science and Technology, Wuhan 430074, China}

\author{Rui Du }\email{rdu@seu.edu.cn}
\affiliation{School of Mathematics, Southeast University, Nanjing 20096, China}

\author{Zhenhua Chai }\email {hustczh@hust.edu.cn}
\affiliation{School of Mathematics and Statistics, Huazhong
University of Science and Technology, Wuhan 430074, China}
\affiliation{Institute of Interdisciplinary Research for Mathematics and Applied Science, Huazhong
University of Science and Technology, Wuhan 430074, China}
\affiliation{Hubei Key Laboratory of Engineering Modeling and
Scientific Computing, Huazhong University of Science and Technology,
Wuhan 430074, China}


\date{\today}

\begin{abstract}
In this paper, based on the idea of direct discrete modeling (DDM) with equilibrium distribution functions (EDFs), we develop a general framework of the mesoscopic numerical method (MesoNM) for macroscopic partial differential equations (PDEs), including but not limited to the nonlinear convection-diffusion equation (NCDE) and the Navier-Stokes equations (NSEs). Unlike the mesoscopic lattice Boltzmann method, this kind of MesoNM is an EDF-based  mesoscopic finite-difference (MesoFD) method, and by taking the moments of the MesoFD scheme, its macroscopic version, called MMFD method,  can be derived directly. 
Both MesoFD scheme and MMFD schemes are multi-level FD methods, MesoFD scheme being mesoscopic, and MMFD scheme being its macroscopic form which has the form of the central FD scheme. They are unified FD schemes for PDEs and can be in implicit or explicit forms as needed.
The macroscopic moment equations (MEs) can be derived from the MesoFD or MMFD scheme through the Taylor expansion method, and the common PDEs can be recovered from the MEs by using the direct Taylor expansion method.
Moreover, the stability of the MMFD scheme is analyzed for linear CDE and liner wave equation with anisotropic diffusion, and the stability conditions of a two-level explicit MMFD scheme, a two-level $\theta$-MMFD scheme (hybrid explicit and implicit MMFD scheme), and a three-level MMFD scheme are obtained, respectively. Finally, we note that some existing lattice Boltzmann (LB) based macroscopic FD models for the NSEs and NCDE are the special cases of present MMFD, which can be considered as a unified framework of FD schemes for PDEs, from this point of view.
\end{abstract}


\maketitle

\section{Introduction}

The lattice Boltzmann (LB) method, as an effective mesoscopic numerical method (MesoNM) based on kinetic theory, has gained great success in the study of complex flows, heat and mass transfer described by the Navier-Stokes equations (NSEs), nonlinear convection-diffusion equations (NCDEs), and their coupled systems \cite{Chen1998,Succi2001,Guo2013,Kruger2017}.
Until now, although there are many different types of
LB models [16], the most popular are the single-relaxation-time LB (SRT-LB) model or lattice Bhatnagar-Gross-Krook (LBGK) model \cite{Qian1990,Chen1991}, the two-relaxation-time LB (TRT-LB) model \cite{Ginzburg2005a,Ginzburg2008,dHumieres2009,Ginzburg2012,Ginzburg2015}, and the general or multple-relaxation-time LB (MRT-LB) model \cite{dHumieres1992,dHumieres2002,Pan2006,Luo2011}. The MRT-LB model has multiple relaxation parameters, and compared to the SRT version, it can achieve more stable and more accurate results by adjusting the free relaxation parameters \cite{dHumieres1992,dHumieres2002,Pan2006,Luo2011}, also the TRT-LB model with only one free relaxation parameter \cite{Ginzburg2005a,Ginzburg2008,dHumieres2009,Ginzburg2012,Ginzburg2015}. In fact, the classic MRT-LB model can be more general by introducing a block-lower-triangular-relaxation matrix \cite{Chai2020}, and rectangular grid \cite{ChaiYuanShi2023}.

It is known that the LB method is based on the evolution of the distribution function (DF) in all lattice velocity directions at each lattice node, and the macroscopic conserved quantities are computed from the moments of  DF. This evolution based on the DF brings us several challenging problems that need to be addressed: (1) a large amount of memory is required; (2) the initial and boundary conditions for macroscopic quantities should be treated carefully; (3) more importantly, the theoretical analysis (e.g., stability analysis) on the LB method is complicated and difficult. To overcome these problems, it is desirable to explore the connection between the LB method and the traditional numerical method. In the past years, some efforts have been made to derive the macroscopic numerical scheme from the LB method.

In the LB community, it is well known that the SRT-LB model with relaxation time $\tau=1$ can be written as a macroscopic FD scheme \cite{Junk1999,Junk2001,Inamuro2002,Sman2006}. In order to overcome the limitation $\tau=1$, the lattice kinetic scheme was proposed in Ref. \cite{Inamuro2002}, and several modified versions have been developed in the literature.
In addition to those derived from the SRT-LB model \cite{Junk1999,Junk2001,Inamuro2002,Sman2006}, up to now there have been many research works on macroscopic FD schemes derived or developed from LB models, e.g., the TRT-LB model \cite{dHumieres2009,Ginzburg2010,Ginzburg2012}, and some other type of macroscopic FD schemes, such as the so-called simplified lattice Boltzmann method \cite{ChenShu2017,ChenShu2018} and its modified models \cite{Gutierrez2021,RosisLiu2021}, and the recursive FD LB models \cite{Vienne2021, HuangChaiShi2022}. Furthermore, for the one-dimensional case, some three- and four-level FD schemes can also be obtained from the LB models \cite{Suga2010,LiZheng2012, ChenChaiShi2023}.
However, it should be noted that based on the recurrence equations with $s^{+}+s^{-}=2$, a class of three-level FD schemes can be derived from TRT-LB model \cite{dHumieres2009,Ginzburg2010,Ginzburg2012}, which contains the FD scheme obtained from the SRT-LB model as its special case (setting $s^{+}=s^{-}=1$), where $s^{+}$ and $s^{-}$ are the relaxation parameters corresponding to the even-order and odd-order moments, respectively.

It can be seen that the above works on the macroscopic FD schemes \cite{FDbook1} derived from the LB models are usually limited to some special cases. Recently, Bellotti et al. \cite{Bellotti2022,Bellotti2023} demonstrated that any LB scheme can be rewritten as a corresponding multi-level FD scheme on the conserved variables, and conducted a rigorous consistency study and the derivation of
the modified equations. However, as stated in Ref. \cite{Bellotti2022}, the multi-level FD scheme usually has the theoretical value due to its complicated multi-level structure and the formalism based on operators, commutative algebra and polynomials ('we stress once more the fact that recasting a lattice Boltzmann scheme as a Finite Difference scheme on the conserved moments is merely useful from the theoretical perspective, being of moderate interest for actual efficient implementations' \cite{Bellotti2022}). Based on the results in Refs. \cite{Bellotti2022,Bellotti2023}, some specific LB models and corresponding multi-level FD schemes with high-order accuracy for one- or two-dimensional simple problems have been designed \cite{ChenChaiShi2023, ChenChaiJCP2024}. According to the previous work \cite{Chai2020}, Chen et al. \cite{ChenChaiPRE2024} presented a general propagation multiple-relaxation-time lattice Boltzmann (GPMRT-LB)
model, and then following the works \cite{Bellotti2022,Bellotti2023}, they obtained the macroscopic FD schemes and modified equations corresponding to the GPMRT-LB model.

Differently from the above methods for deriving the FD schemes from the existing LB models, in this work, we use the equilibrium distribution function (EDF) to represent the distribution function (DF), and propose a novel class of mesoscopic numerical methods (MesoNMs), called the mesoscopic finite-difference (MesoFD) schemes, for partial differential equations (PDEs) with the source terms. Then the MesoFD based macroscopic finite-difference (MMFD) schemes can be directly obtained by taking the moments of the MesoFD scheme. Both of them are of multi-level implicit or explicit structure. The MMFD scheme can be viewed as a unified framework of FD schemes for PDEs, and the commonly used FD schemes derived from LB models and traditional FD schemes are its special cases.

The rest of this paper is organized as follows. In Sec. II, the MesoFD method and the corresponding MMFD scheme are presented. The direct Taylor expansion is adopted to obtain the corresponding moment equations and PDEs of the present MesoFD and MMFD methods, and we present some special cases of the MesoFD and MMFD methods in Sec. III, and followed by some numerical tests in Sec. IV. Finally, some conclusions are summarized in Sec. V.

\section{Mesoscopic finite difference method and macroscopic finite difference scheme}

The MesoFD method comes from the idea of representing the distribution function (DF) or nonequilibrium distribution function (NEDF) with equilibrium distribution function (EDF), which is an extension of the previous work \cite{HuangChaiShi2022}.

The evolution equation of MesoFD method with the general rectangular D$d$Q$b$ (rD$d$Q$b$) lattice \cite{ChaiYuanShi2023} has the following multi-level form:
\begin{equation}\label{eq:2-1}
f_j(\mathbf{x},t)=\sum_{k,q}\big[a_{kq}f_j^{eq[k,q]}+\Delta t b_{kq}F_j^{[k,q]}\big],
\end{equation}
with
\begin{equation}\label{eq:2-1-0}
f_j^{eq[k,q]}=f_j^{eq}(\mathbf{x}+k\mathbf{c}_j \Delta t,t+q \Delta t), \ F_j^{[k,q]}=F_j(\mathbf{x}+k\mathbf{c}_j \Delta t,t+q \Delta t), \ j=0, 1, \cdots, b-1,
\end{equation}
where $\mathbf{c}_j$ is the discrete velocity in rD$d$Q$b$ lattice model, and the details can be found in Ref. \cite{ChaiYuanShi2023}, while for convenience, some commonly used lattice models are given in Appendix A. $f_j(\mathbf{x}, t)$ is the distribution function at position $\mathbf{x}$ in $d$-dimensional space and time $t$ along the velocity $\mathbf{c}_j$, and $f_j^{eq}(\mathbf{x}, t)$ is the equilibrium distribution function. Let $f_j^{ne}(\mathbf{x}, t)=f_j(\mathbf{x}, t)-f_j^{eq}(\mathbf{x}, t)$, which is the nonequilibrium distribution function (NEDF). $F_j(\mathbf{x}, t)$ is the distribution function of the source or forcing term. The parameters $a_{kq}$ and $b_{kq}$ need to be determined, where $k$ and $q$ are integers in a given finite set $\{-m_q\leq k\leq m_q, -n\leq q\leq 0\}(n\geq 1,m_q\geq 1 $ for all $q$ ). $\Delta t$ is the time step.
In the evolution equation (\ref{eq:2-1}), the key elements, $\mathbf{c}_j, f_{j}^{eq}$, $F_j$, and the parameters $a_{kq}$ and $b_{kq}$, must be given properly.

From Eq. (\ref{eq:2-1}) we can also obtain the evolution equations of required moments. Similar to the LB method, the unknown macroscopic conserved variable(s), for instance, the zeroth moment $\phi(\mathbf{x},t)$ for scalar PDE, such as NCDE, or the zeroth moment $\rho(\mathbf{x},t)$ and the first moment $(\rho\mathbf{u})(\mathbf{x},t)$ for NSEs, can be computed by

\begin{subequations}\label{eq:2-1-1}
\begin{equation}
\phi(\mathbf{x}, t)=\sum_j f_j(\mathbf{x}, t),
\end{equation}
\begin{equation}
\rho(\mathbf{x}, t)=\sum_j f_j(\mathbf{x}, t),
\quad(\rho\mathbf{u})(\mathbf{x}, t)=\sum_j \mathbf{c}_j f_j(\mathbf{x}, t).
\end{equation}
\end{subequations}

In this case, the corresponding MMFD scheme can be obtained from Eq. (1),
\begin{equation}\label{eq:2-2}
\textbf{NCDE:}\ \   \phi(\mathbf{x},t)=\sum_j\sum_{k,q}\big[a_{kq}f_j^{eq[k,q]}+\Delta t b_{kq}F_j^{[k,q]}\big];  \\
\end{equation}
\begin{subequations}\label{eq:2-3}
\begin{equation}
\textbf{NSEs:}\ \   \rho(\mathbf{x},t)=\sum_j\sum_{k,q}\big[a_{kq}f_j^{eq[k,q]}+\Delta t b_{kq}F_j^{[k,q]}\big],  \ \
\end{equation}
\begin{equation}
\ \ \ \ \ \ \ \ \   (\rho \mathbf{u})(\mathbf{x},t)=\sum_j\sum_{k,q}\mathbf{c}_j\big[a_{kq} f_j^{eq[k,q]}+\Delta t b_{kq} F_j^{[k,q]}\big].  \ \
\end{equation}
\end{subequations}

Furthermore, through using a \textit{divide-and-conquer} strategy, we can extend the single-distribution-function scheme Eq. (\ref{eq:2-3}) to a double-distribution-function one:
\begin{subequations}\label{eq:2-4}
\begin{equation}
\rho(\mathbf{x},t)=\sum_j\sum_{k,q}\big[a_{kq}f_j^{eq[k,q]}+\Delta t b_{kq}F_j^{[k,q]}\big],  \ \
\end{equation}
\begin{equation}
(\rho\mathbf{u})(\mathbf{x},t)=\sum_j\sum_{k,q}\mathbf{c}_j\big[\bar{a}_{kq} f_j^{eq[k,q]}+\Delta t \bar{b}_{kq} F_j^{[k,q]}\big],  \ \
\end{equation}
\end{subequations}
in which two different sets of parameters, $\{a_{kq},b_{kq}\}$ and $\{\bar{a}_{kq},\bar{b}_{kq}\}$, are used.

\textbf{Remark 1} We can replace $f_j$ and $a_{00}$ in Eq. (\ref{eq:2-1}) by $f_j^{ne}$ and $a_{00}-1$, respectively, then we obtain the equivalent form of Eq. (\ref{eq:2-1}). Therefore, if the summation is only made for the case $k=q$, one can easily find that the scheme in Ref. \cite{HuangChaiShi2022} is just a special case of the MesoFD scheme (\ref{eq:2-1}).  Unless otherwise stated, the following discussion is only limited to Eq. (\ref{eq:2-1}).

\subsection{The central difference representation of MesoFD method}

It is interesting that both MesoFD and MMFD schemes can be represented by the central difference format. First, we define the following central difference operators through introducing a function $g(\mathbf{x})$, and two integers $j$ and $k$, 
\begin{eqnarray}\label{eq:CentralFD}
&\Delta_{j,k}^{(1)}g(\mathbf{x})=g(\mathbf{x}+k\mathbf{c}_j\Delta t)-g(\mathbf{x}-k\mathbf{c}_j\Delta t),\nonumber\\
&\Delta_{j,k}^{(2)}g(\mathbf{x})=g(\mathbf{x}+k\mathbf{c}_j\Delta t)-2g(\mathbf{x})+g(\mathbf{x}-k\mathbf{c}_j\Delta t),\nonumber\\
&\Delta_{j}^{(1)}g(\mathbf{x})=\Delta_{j,1}^{(1)}g(\mathbf{x})=g(\mathbf{x}+\mathbf{c}_j\Delta t)-g(\mathbf{x}-\mathbf{c}_j\Delta t), \nonumber\\
&\Delta_{j}^{(2)}g(\mathbf{x})=\Delta_{j,1}^{(2)}g(\mathbf{x})=g(\mathbf{x}+\mathbf{c}_j\Delta t)-2g(\mathbf{x})+g(\mathbf{x}-\mathbf{c}_j\Delta t).
\end{eqnarray}
It is easy to obtain
\begin{eqnarray}\label{eq:CentralFD-0}
\Delta_{j,k}^{(1)}=-\Delta_{j,-k}^{(1)}=-\Delta_{\bar{j},k}^{(1)}, \quad \Delta_{j,k}^{(2)}=\Delta_{j,-k}^{(2)}=\Delta_{\bar{j},k}^{(2)},\nonumber\\
\Delta_{0,k}^{(1)}=\Delta_{0,k}^{(2)}=\Delta_{j,0}^{(1)}=\Delta_{j,0}^{(2)}=0,
\end{eqnarray}
where $\bar{j}$ is the opposite direction of $j$,  $\mathbf{c}_{\bar{j}}=-\mathbf{c}_j$, $\mathbf{c}_0=\mathbf{0}$.

Using the above operators, one can obtain the central difference representation of MesoFD scheme (\ref{eq:2-1}):
\begin{eqnarray}
f_j(\mathbf{x},t)=\sum_{k,q}a_{kq}\big[1+\frac{1}{2}(\Delta_{j,k}^{(2)}+\Delta_{j,k}^{(1)})\big] f_j^{eq}(\mathbf{x},t+q\Delta t)
 +\Delta t \sum_{k,q}b_{kq}\big[1+\frac{1}{2}(\Delta_{j,k}^{(2)}+\Delta_{j,k}^{(1)})\big] F_j(\mathbf{x},t+q\Delta t),\nonumber\\
 \label{eq:2-1-FDa}
\\
 f_{\bar{j}}(\mathbf{x},t)=\sum_{k,q}a_{kq}\big[1+\frac{1}{2}(\Delta_{j,k}^{(2)}-\Delta_{j,k}^{(1)})\big] f_{\bar{j}}^{eq}(\mathbf{x},t+q\Delta t)
 +\Delta t \sum_{k,q}b_{kq}\big[1+\frac{1}{2}(\Delta_{j,k}^{(2)}-\Delta_{j,k}^{(1)})\big] F_{\bar{j}}(\mathbf{x},t+q\Delta t).\nonumber\\
  \label{eq:2-1-FDb}
\end{eqnarray}


For a given distribution function $g_j$, let
\begin{equation}\label{eq:2-5}
g_j^{+}=\frac{1}{2}(g_j + g_{\bar{j}}), \ g_j^{-}=\frac{1}{2}(g_j - g_{\bar{j}}),
\end{equation}
then we have
\begin{equation}\label{eq:2-6}
g_j^{+}=g_{\bar{j}}^{+}, \ g_j^{-}=-g_{\bar{j}}^{-},\ g_j=g_j^{+}+g_j^{-}, \ 
g_{\bar{j}}=g_j^{+}-g_j^{-}.
\end{equation}

With the help of Eqs. (\ref{eq:CentralFD-0}), (\ref{eq:2-5}) and (\ref{eq:2-6}), the summation and the difference of Eqs. (\ref{eq:2-1-FDa}) and (\ref{eq:2-1-FDb})  become
\begin{subequations}\label{eq:2-1-FD1}
\begin{eqnarray}
 f_j^{+}(\mathbf{x},t)&=&\sum_{k,q}a_{kq}\big[f_j^{eq,+}+\frac{1}{2}(\Delta_{j,k}^{(2)}f_j^{eq,+}+\Delta_{j,k}^{(1)}f_j^{eq,-})\big] (\mathbf{x},t+q\Delta t)
 +\Delta t \sum_{k,q}b_{kq}\big[F_j^{+}+\frac{1}{2}(\Delta_{j,k}^{(2)}F_j^{+}+\Delta_{j,k}^{(1)}F_j^{-})\big] (\mathbf{x},t+q\Delta t)\nonumber\\
& =&\sum_{k,q}\big[a_{kq}f_j^{eq,+}+ \Delta t b_{kq}F_j^{+} \big](\mathbf{x},t+q\Delta t) 
+\frac{1}{2}\sum_q \sum_{k=1}^{m_q}\big[ a_{kq}^{+}\Delta_{j,k}^{(2)}f_j^{eq,+}+a_{kq}^{-}\Delta_{j,k}^{(1)}f_j^{eq,-})\big] (\mathbf{x},t+q\Delta t)\nonumber\\
& & + \frac{\Delta t}{2}\sum_q \sum_{k=1}^{m_q}\big[ b_{kq}^{+}\Delta_{j,k}^{(2)}F_j^{+}+b_{kq}^{-}\Delta_{j,k}^{(1)}F_j^{,-})\big] (\mathbf{x},t+q\Delta t),
\end{eqnarray}
\begin{eqnarray}
 f_j^{-}(\mathbf{x},t)&=&\sum_{k,q}a_{kq}\big[f_j^{eq,-}+\frac{1}{2}(\Delta_{j,k}^{(2)}f_j^{eq,-}+\Delta_{j,k}^{(1)}f_j^{eq,+})\big] (\mathbf{x},t+q\Delta t)
 +\Delta t \sum_{k,q}b_{kq}\big[F_j^{-}+\frac{1}{2}(\Delta_{j,k}^{(2)}F_j^{-}+\Delta_{j,k}^{(1)}F_j^{+})\big] (\mathbf{x},t+q\Delta t),\nonumber\\
& =&\sum_{k,q}\big[a_{kq}f_j^{eq,-}+ \Delta t b_{kq}F_j^{-} \big](\mathbf{x},t+q\Delta t) 
+\frac{1}{2}\sum_q \sum_{k=1}^{m_q}\big[ a_{kq}^{+}\Delta_{j,k}^{(2)}f_j^{eq,-}+a_{kq}^{-}\Delta_{j,k}^{(1)}f_j^{eq,+})\big] (\mathbf{x},t+q\Delta t)\nonumber\\
& & + \frac{\Delta t}{2}\sum_q \sum_{k=1}^{m_q}\big[ b_{kq}^{+}\Delta_{j,k}^{(2)}F_j^{-}+b_{kq}^{-}\Delta_{j,k}^{(1)}F_j^{,+})\big] (\mathbf{x},t+q\Delta t),
\end{eqnarray}
\end{subequations}
where
\begin{equation}\label{eq:2-1-FD1-0}
a_{kq}^{+}=a_{kq}+a_{-kq}, ~a_{kq}^{-}=a_{kq}-a_{-kq}, ~b_{kq}^{+}=b_{kq}+b_{-kq}, ~b_{kq}^{-}=b_{kq}-b_{-kq}.
\end{equation}

Based on the following equations, 
\begin{eqnarray}
\sum_j h_j^{+} = \sum_j h_j, \sum_j h_j^{-}=0, ~\sum_j \mathbf{c}_j h_j^{+} = \mathbf{0}, ~\sum_j \mathbf{c}_j h_j^{-}=\sum_j \mathbf{c}_j h_j, ~h_j \in \{f_j, f_j^{eq}, F_j\},
\end{eqnarray}
we can take the zeroth moment of $f_j^{+}$ and the first moment of $f_j^{-}$ in Eq. (\ref{eq:2-1-FD1}), respectively, and get the MMFD schemes for the NCDE and NSEs, which are equivalent to those in Eqs. (\ref{eq:2-2}) and (\ref{eq:2-4})
\begin{eqnarray}\label{eq:2-2-1}
    \phi(\mathbf{x},t)&=&\sum_{k,q}\big[a_{kq} \phi+\Delta t b_{kq}S\big](\mathbf{x},t+q\Delta t)
    + \frac{1}{2}\sum_q\sum_j\sum_{k=1}^{m_q}\big[a_{kq}^{+}\Delta_{j,k}^{(2)}f_j^{eq,+}+\Delta t b_{kq}^{+}\Delta_{j,k}^{(2)}F_j^{+}\big](\mathbf{x},t+q\Delta t)\nonumber\\
& &+  \frac{1}{2}\sum_q\sum_j\sum_{k=1}^{m_q}\big[a_{kq}^{-}\Delta_{j,k}^{(1)}f_j^{eq,-}+\Delta t b_{kq}^{-}\Delta_{j,k}^{(1)}F_j^{-}\big](\mathbf{x},t+q\Delta t);
\end{eqnarray}
\begin{subequations}\label{eq:2-4-1}
\begin{eqnarray}
    \rho(\mathbf{x},t)&=&\sum_{k,q}\big[a_{kq} \rho+\Delta t b_{kq}S\big](\mathbf{x},t+q\Delta t)
   + \frac{1}{2}\sum_q\sum_j\sum_{k=1}^{m_q}\big[a_{kq}^{+}\Delta_{j,k}^{(2)}f_j^{eq,+}+\Delta t b_{kq}^{+}\Delta_{j,k}^{(2)}F_j^{+}\big](\mathbf{x},t+q\Delta t)\nonumber\\
& &+  \frac{1}{2}\sum_q\sum_j\sum_{k=1}^{m_q}\big[a_{kq}^{-}\Delta_{j,k}^{(1)}f_j^{eq,-}+\Delta t b_{kq}^{-}\Delta_{j,k}^{(1)}F_j^{-}\big](\mathbf{x},t+q\Delta t),
\end{eqnarray}
\begin{eqnarray}
    (\rho \mathbf{u})(\mathbf{x},t)&=&\sum_{k,q}\big[\bar{a}_{kq} (\rho \mathbf{u})+\Delta t \bar{b}_{kq}\mathbf{F}\big](\mathbf{x},t+q\Delta t)
    + \frac{1}{2}\sum_q\sum_j\sum_{k=1}^{m_q}\mathbf{c}_j\big[\bar{a}_{kq}^{+}\Delta_{j,k}^{(2)}f_j^{eq,-}+\Delta t \bar{b}_{kq}^{-}\Delta_{j,k}^{(2)}F_j^{-}\big](\mathbf{x},t+q\Delta t)\nonumber\\
& &+  \frac{1}{2}\sum_q\sum_j\sum_{k=1}^{m_q}\mathbf{c}_j\big[\bar{a}_{kq}^{-}\Delta_{j,k}^{(1)}f_j^{eq,+}+\Delta t \bar{b}_{kq}^{-}\Delta_{j,k}^{(1)}F_j^{+}\big](\mathbf{x},t+q\Delta t),
\end{eqnarray}
\end{subequations}
where $S,\mathbf{F}$ in Eqs. (\ref{eq:2-2-1}) and (\ref{eq:2-4-1}) and the parameters in Eq. (\ref{eq:2-4-1}b) are defined as $S=\sum_j F_j=\sum_j F_j^{+},~\mathbf{F}=\sum_j \mathbf{c}_j F_j=\sum_j \mathbf{c}_j F_j^{-},~$ $\bar{a}_{kq}^{+}=\bar{a}_{kq}+\bar{a}_{-kq},~\bar{a}_{kq}^{-}=\bar{a}_{kq}-\bar{a}_{-kq},~\bar{b}_{kq}^{+}=\bar{b}_{kq}+\bar{b}_{-kq},~\bar{b}_{kq}^{-}=\bar{b}_{kq}-\bar{b}_{-kq}$, $\bar{a}_{kq} and \bar{b}_{kq}$  are used as those in Eq. (\ref{eq:2-4}). If needed, one can take the higher-order moments in Eq. (\ref{eq:2-1-FD1}) to get the evolution equations of these higher-order moments.

It can be seen from Eqs. (\ref{eq:2-2-1}) and (\ref{eq:2-4-1}) that the MMFD scheme consists of central difference schemes at multiple time levels. Based on the MesoFD scheme (Eq. (\ref{eq:2-1})) and the MMFD scheme (Eq. (\ref{eq:2-2-1}) or Eq. (\ref{eq:2-4-1})), we can derive the macroscopic PDEs from mesoscopic and macroscopic ways, respectively.

It is worth noting that both Eq. (\ref{eq:2-2-1}) and Eq. (\ref{eq:2-4-1}) have the seperated form of spatio-temporal difference. Through using a \textit{divide-and-conquer} strategy again, we can extend the one-set-parameter scheme Eq. (\ref{eq:2-2-1})  to a multiple-set-parameter one,
\begin{eqnarray}\label{eq:2-2-1b}
    \phi(\mathbf{x},t)&=&\sum_{k,q}\big[a_{kq} \phi+\Delta t b_{kq}S\big](\mathbf{x},t+q\Delta t)
    + \frac{1}{2}\sum_q\sum_j\sum_{k=1}^{m_q}\big[\bar{a}_{kq}^{+}\Delta_{j,k}^{(2)}f_j^{eq,+}+\Delta t \bar{b}_{kq}^{+}\Delta_{j,k}^{(2)}F_j^{+}\big](\mathbf{x},t+q\Delta t)\nonumber\\
& &+  \frac{1}{2}\sum_q\sum_j\sum_{k=1}^{m_q}\big[\hat{a}_{kq}^{-}\Delta_{j,k}^{(1)}f_j^{eq,-}+\Delta t \hat{b}_{kq}^{-}\Delta_{j,k}^{(1)}F_j^{-}\big](\mathbf{x},t+q\Delta t),
\end{eqnarray}
 and correspondingly, the MesoFD scheme (1) has the following form,
\begin{eqnarray}\label{eq:2-1b}
f_j(\mathbf{x},t)=\sum_{k,q}\big[\bar{a}_{kq}f_j^{eq[k,q],+}+\hat{a}_{kq}f_j^{eq[k,q],-}+\Delta t (\bar{b}_{kq}F_j^{[k,q],+}+\hat{b}_{kq}F_j^{[k,q],-})\big],
\end{eqnarray}
where $\{a_{kq},b_{kq}\},\{\bar{a}_{kq},\bar{b}_{kq}\}$ and $\{\hat{a}_{kq},\hat{b}_{kq}\}$ can be different. Similarly, one can give the multiple-set-parameter version of scheme Eq. (\ref{eq:2-4-1}). It should be noted that Eq. (\ref{eq:2-1b}) is a TRT-version of the MesoFD scheme (1), and one can also give a MRT-version of the MesoFD scheme (1) by assigning different coefficients to each term of $f_j^{eq}$. For simplicity, however, in the following, we will not discuss the TRT- and MRT- versions.

\subsection{The two-level central MMFD scheme}
As we know, the two-level FD scheme is the most widely used one for its efficiency. In this part, we only consider the two-level central MMFD scheme for NCDE, and the two-level central MMFD scheme for NSEs can be given in a similar way. Taking $q=-1$ in Eq. (\ref{eq:2-2-1}), we obtain the explicit MMFD (eMMFD) scheme, 
\begin{eqnarray}\label{eq:2-2-2}
    \phi(\mathbf{x},t)&=&\sum_{k}\big[a_{k,-1} \phi+\Delta t b_{k,-1}S\big](\mathbf{x},t-\Delta t)
    + \frac{1}{2}\sum_j\sum_{k=1}^{m_{-1}}\big[a_{k,-1}^{+}\Delta_{j,k}^{(2)}f_j^{eq,+}+\Delta t b_{k,-1}^{+}\Delta_{j,k}^{(2)}F_j^{+}\big](\mathbf{x},t-\Delta t)\nonumber\\
& & +  \frac{1}{2}\sum_j\sum_{k=1}^{m_{-1}}\big[a_{k,-1}^{-}\Delta_{j,k}^{(1)}f_j^{eq,-}+\Delta t b_{k,-1}^{-}\Delta_{j,k}^{(1)}F_j^{-}\big](\mathbf{x},t-\Delta t),
\end{eqnarray}
and taking $q \in \{0,-1\}$ in Eq. (\ref{eq:2-2-1}), the implicit MMFD (iMMFD) scheme  can be obtained as follows,
\begin{eqnarray}\label{eq:2-2-3}
    \phi(\mathbf{x},t)&=&\sum_{k}\big[a_{k,-1} \phi+\Delta t b_{k,-1}S\big](\mathbf{x},t-\Delta t) + \frac{1}{2}\sum_j\sum_{k=1}^{m_0}\big[a_{k,0}^{+}\Delta_{j,k}^{(2)}f_j^{eq,+}+\Delta t b_{k,0}^{+}\Delta_{j,k}^{(2)}F_j^{+}\big](\mathbf{x},t)\nonumber\\
& &+  \frac{1}{2}\sum_j\sum_{k=1}^{m_0}\big[a_{k,0}^{-}\Delta_{j,k}^{(1)}f_j^{eq,-}+\Delta t b_{k,0}^{-}\Delta_{j,k}^{(1)}F_j^{-}\big](\mathbf{x},t),
\end{eqnarray}
The eMMFD and iMMFD schemes can be combined to give a more general two-level $\theta$-scheme, 
\begin{eqnarray}\label{eq:GCN}
   & & \phi(\mathbf{x},t)=\sum_{k}\big[a_{k} \phi+\Delta t b_{k}S\big](\mathbf{x},t-\Delta t)\nonumber\\
  & & + \frac{1}{2}\sum_j\sum_{k=1}^{m_{0}}\big[a_{k}^{+}\Delta_{j,k}^{(2)}(\theta f_j^{eq,+}(\mathbf{x},t-\Delta t)+\theta_1 f_j^{eq,+}(\mathbf{x},t))+\Delta t b_{k}^{+}\Delta_{j,k}^{(2)}(\theta F_j^{+}(\mathbf{x},t-\Delta t)+\theta_1 F_j^{+}(\mathbf{x},t))\big]\nonumber\\
& &+  \frac{1}{2}\sum_j\sum_{k=1}^{m_{0}}\big[a_{k}^{-}\Delta_{j,k}^{(1)}(\theta f_j^{eq,-}(\mathbf{x},t-\Delta t)+\theta_1 f_j^{eq,-}(\mathbf{x},t))+\Delta t b_{k}^{-}\Delta_{j,k}^{(1)}(\theta F_j^{-}(\mathbf{x},t-\Delta t)+\theta_1 F_j^{-}(\mathbf{x},t))\big],
\end{eqnarray}
where $\theta \in [0,1],~\theta_1=1-\theta, ~a_{k}^{+}=a_{k}+a_{-k},~a_{k}^{-}=a_{k}-a_{-k},$ and $b_{k}^{+}=b_{k}+b_{-k},b_{k}^{-}=b_{k}-b_{-k}$.

\section{The moment equations of MesoFD method: Direct Taylor expansion}
\label{printlayout}

There are four basic analysis methods that can be used to recover the macroscopic NSEs and NCDE from the LB models, i.e., the Chapman-Enskog (C-E) analysis \cite{Chapman1970}, the Maxwell iteration (MI) method \cite{Yong2016}, the direct Taylor expansion (DTE) method \cite{Chai2020} and the recurrence equations (RE) method \cite{Ginzburg2012,dHumieres2009,Ginzburg2015}. In Ref.  \cite{Chai2020}, it is shown that they can give the same equations at the second-order of expansion parameters, while the DTE method is much simpler. The MI method and RE method may not be suitable for analyzing the MesoFD scheme. In what follows, the DTE method is used to analyze the MesoFD scheme at the acoustic scaling.

Applying the Taylor expansion to Eq. (\ref{eq:2-1}),  one can get
\begin{equation}\label{eq:4-1}
f_j= \sum_{l=0}^{N}\frac{\Delta t^l}{l!}  \sum_{k,q}\big[a_{kq}D_{k,q}^l f_j^{eq}+\Delta t b_{kq}D_{k,q}^l F_j\big] +O(\Delta t^{N+1}), ~N\geq 1,
\end{equation}
where $D_{k,q}=q\partial_t + k \nabla\cdot\mathbf{c}_j$.

Let
\begin{equation}\label{eq:4-2}
A_{lm}=\frac{1}{l!}C_l^m  \sum_{k,q}a_{kq}q^{l-m}k^m,  ~B_{lm}=\frac{1}{l!}C_l^m  \sum_{k,q}b_{kq}q^{l-m}k^m,
\end{equation}
where $C_l^m=\frac{l!}{m!(l-m)!}$. From Eq. (\ref{eq:4-1}), we can derive the following equations,
\begin{equation}\label{eq:4-1-1}
f_j= \sum_{l=0}^{N}\Delta t^l\sum_{m=0}^{l}\big[A_{lm}\partial_t^{l-m}(\nabla\cdot \mathbf{c}_j)^m f_j^{eq}+\Delta t B_{lm}\partial_t^{l-m}(\nabla\cdot \mathbf{c}_j)^m F_j\big] +O(\Delta t^{N+1}), ~N\geq 1,
\end{equation}
where the following formula is used
\begin{equation}\label{eq:4-2-1}
D_{k,q}^l=(q\partial_t + k \nabla\cdot\mathbf{c}_j)^l=\sum_{m=0}^{l} C_l^m  (q\partial_t)^{l-m}(k\nabla\cdot \mathbf{c}_j)^m.
\end{equation}

Taking $N=$1 and 2 in Eq. (\ref{eq:4-1-1}), one can obtain
\begin{subequations}\label{eq:4-3}
\begin{equation}
f_j= A_{00}f_j^{eq}+\Delta t(A_{10}\partial_t+A_{11}\nabla\cdot \mathbf{c}_j) f_j^{eq}+\Delta t B_{00} F_j +O(\Delta t^{2}),
\end{equation}
and
\begin{eqnarray}
f_j&=& A_{00}f_j^{eq}+\Delta t(A_{10}\partial_t+A_{11}\nabla\cdot \mathbf{c}_j) f_j^{eq}+\Delta t^2\big[A_{20}\partial_t^2+A_{21}\partial_t \nabla\cdot \mathbf{c}_j+A_{22}(\nabla\cdot \mathbf{c}_j)^2\big] f_j^{eq}\nonumber\\
& &\Delta t B_{00} F_j +\Delta t^2 (B_{10}\partial_t+B_{11}\nabla\cdot \mathbf{c}_j) F_j+O(\Delta t^{3}),
\end{eqnarray}
\end{subequations}
respectively.

\subsection{The derivation of the moment equation for PDE using DTE method}

Based on Eq. (\ref{eq:4-1-1}), we can first derive the moment equation for PDE from the MesoFD method (\ref{eq:2-1}). To simplify the following analysis, we denote the moments of $f_j$, $f_j^{eq}$, and $F_j$ as
\begin{subequations}\label{eq:M-NCDE}
\begin{equation}
M_0=\sum_j f_j^{eq}=\sum_j f_j,\ \  \mathbf{M}_1=\sum_j \mathbf{c}_j f_j^{eq}, \ \  \mathbf{M}_k=\sum_j \mathbf{c}_j \cdot\cdot\cdot\mathbf{c}_j f_j^{eq}, ~k>1,
\end{equation}
\begin{equation}
M_{0F}=\sum_j F_j,\ \ \mathbf{M}_{1F}=\sum_j \mathbf{c}_j F_j, \ \  \mathbf{M}_{kF}=\sum_j \mathbf{c}_j \cdot\cdot\cdot\mathbf{c}_j F_j, ~k>1,
\end{equation}
\end{subequations}
where $\mathbf{M}_k$ and $\mathbf{M}_{kF}$ ($k\geq 0$) are the $k$-th moments of $f_j^{eq}$ and $F_j$, respectively.
Additionally, from Eq. (\ref{eq:M-NCDE}a) one can derive the following moment of nonequilibrium distribution function,
\begin{equation}\label{eq:M-NCDE-1}
M_0^{ne}=\sum_j f_j^{ne}=\sum_j (f_j-f_j^{eq})=0.
\end{equation}

Summing Eq. (\ref{eq:4-1-1}) and adopting Eq. (\ref{eq:M-NCDE}), one can obtain
\begin{eqnarray}\label{eq:ME-NCDE}
M_0&=&\sum_{l=0}^{N}\Delta t^l\sum_{m=0}^{l}\big[A_{lm}\partial_t^{l-m}( \nabla^m \otimes\mathbf{M}_m)+\Delta t B_{lm}\partial_t^{l-m}(\nabla^m \otimes\mathbf{M}_{mF})\big] +O(\Delta t^{N+1})\nonumber\\
&=&A_{00}M_0+\sum_{l=1}^{N}\Delta t^l\sum_{m=0}^{l}\big[A_{lm}\partial_t^{l-m} (\nabla^m \otimes\mathbf{M}_m)\big]+ \sum_{l=0}^{N-1}\Delta t^{l+1}\sum_{m=0}^{l}\big[B_{lm}\partial_t^{l-m} (\nabla^m \otimes\mathbf{M}_{mF})\big]+O(\Delta t^{N+1}), ~N\geq 1, \nonumber \\
\end{eqnarray}
where
\begin{equation}\label{eq:ME-NCDE-Grad}
\nabla^m\otimes \mathbf{M}_m=\sum_j (\nabla\cdot\mathbf{c}_j)^m  f_j^{eq}=(\nabla\cdot\cdot\cdot\nabla)\otimes\sum_j \mathbf{c}_j\cdot\cdot\cdot\mathbf{c}_j  f_j^{eq}.
\end{equation}

Taking $A_{00}=1$ in Eq. (\ref{eq:ME-NCDE}), we can obtain a general moment equation for PDE of order $N$,
\begin{eqnarray}\label{eq:ME-NCDE-0}
\sum_{l=1}^{N}\Delta t^{l-1}\sum_{m=0}^{l}\big[A_{lm}\partial_t^{l-m} (\nabla^m \otimes\mathbf{M}_m)\big]+ \sum_{l=0}^{N-1}\Delta t^{l}\sum_{m=0}^{l}\big[B_{lm}\partial_t^{l-m} (\nabla^m \otimes\mathbf{M}_{mF})\big]=O(\Delta t^{N}), ~N\geq 1,
\end{eqnarray}
which yields for $N=1,2$, and $3$,
\begin{subequations}\label{eq:ME-NCDE-1}
\begin{eqnarray}
A_{10}\partial_t M_0+A_{11}\nabla\cdot \mathbf{M}_1 + B_{00}M_{0F}=O(\Delta t),
\end{eqnarray}
\begin{eqnarray}
A_{10}\partial_t M_0+A_{11}\nabla\cdot \mathbf{M}_1 +B_{00}M_{0F} +\Delta t\big[A_{20}\partial_t^2 M_0+A_{21}\partial_t \nabla\cdot \mathbf{M}_1+A_{22}\nabla\nabla : \mathbf{M}_2\big]\nonumber\\
+\Delta t(B_{10}\partial_t M_{0F}+B_{11}\nabla\cdot \mathbf{M}_{1F} )=O(\Delta t^2),
\end{eqnarray}
\begin{eqnarray}
& &A_{10}\partial_t M_0+A_{11}\nabla\cdot \mathbf{M}_1 +B_{00}M_{0F} +\Delta t\big[A_{20}\partial_t^2 M_0+A_{21}\partial_t \nabla\cdot \mathbf{M}_1+A_{22}\nabla\nabla : \mathbf{M}_2\big]\nonumber\\
& &+\Delta t(B_{10}\partial_t M_{0F}+B_{11}\nabla\cdot \mathbf{M}_{1F} )+\Delta t^2 \big[A_{30}\partial_t^3 M_0+A_{31}\partial_t^2 \nabla\cdot \mathbf{M}_1+A_{32}\partial_t\nabla\nabla : \mathbf{M}_2+A_{33}\nabla^3 \otimes \mathbf{M}_3\big]\nonumber\\
& &+\Delta t^2(B_{20}\partial_t^2 M_{0F}+B_{21}\partial_t\nabla\cdot \mathbf{M}_{1F} +B_{22}\nabla\nabla : \mathbf{M}_{2F})=O(\Delta t^3),
\end{eqnarray}
\end{subequations}
where
\begin{eqnarray}\label{eq:4-2-0}
&A_{00}=\sum_{k,q} a_{kq}=1, ~A_{10}=\sum_{k,q}q \cdot a_{kq}, ~A_{11}=\sum_{k,q}k\cdot a_{kq}, \nonumber\\
&A_{20}=\frac{1}{2}\sum_{k,q}q^2\cdot a_{kq}, ~A_{21}=\sum_{k,q}q\cdot k\cdot a_{kq},  ~A_{22}=\frac{1}{2}\sum_{k,q}k^2\cdot a_{kq},\nonumber\\
&A_{30}=\frac{1}{6}\sum_{k,q}q^3\cdot a_{kq}, ~A_{31}=\frac{1}{2}\sum_{k,q}q^2\cdot k\cdot a_{kq}, ~ A_{32}=\frac{1}{2}\sum_{k,q}q \cdot k^2\cdot a_{kq}, ~A_{33}=\frac{1}{6}\sum_{k,q}k^3\cdot a_{kq}, \nonumber\\
&B_{00}=\sum_{k,q} b_{kq}, ~B_{10}=\sum_{k,q}q \cdot b_{kq}, ~B_{11}=\sum_{k,q}k\cdot b_{kq},\nonumber\\
&B_{20}=\frac{1}{2}\sum_{k,q}q^2\cdot b_{kq}, ~B_{21}=\sum_{k,q}q\cdot k\cdot b_{kq},  ~B_{22}=\frac{1}{2}\sum_{k,q}k^2\cdot b_{kq}.
\end{eqnarray}

We note that Eq. (\ref{eq:ME-NCDE-1}c) is the moment equation corresponding to a third-order PDE, while Eq. (\ref{eq:ME-NCDE-1}a) and Eq. (\ref{eq:ME-NCDE-1}b) are the first- and second-order approximation, and is also a moment equation for second-order PDE. Some popular second-order PDEs, such as NCDE, convection equation (CE), Poisson equation, and even wave equation (WE) can be derived from Eq. (\ref{eq:ME-NCDE-1}b) or Eq. (\ref{eq:ME-NCDE-1}c) when $A_{lm}, B_{lm}, \mathbf{\mathbf{M}}_k$, and $\mathbf{M}_{kF}$  are properly given.

(1) NCDE

When $M_0$ and $\mathbf{M}_1$ are not constants, taking $A_{10}=A_{11}$ ($\neq 0$) in Eq. (\ref{eq:ME-NCDE-1}), we have
\begin{subequations}\label{eq:ME-NCDE-2}
\begin{eqnarray}
\partial_t M_0+\nabla\cdot \mathbf{M}_1=-\frac{B_{00}}{A_{10}}M_{0F}+O(\Delta t),
\end{eqnarray}
\begin{eqnarray}
\partial_t M_0+\nabla\cdot \mathbf{M}_1+ \frac{\Delta t}{A_{10}}\big[A_{20}\partial_t^2 M_0+A_{21}\partial_t \nabla\cdot \mathbf{M}_1+A_{22}\nabla\nabla : \mathbf{M}_2\big]\nonumber\\
=-\frac{B_{00}}{A_{10}}M_{0F}-\frac{\Delta t}{A_{10}}(B_{10}\partial_t M_{0F}+B_{11}\nabla\cdot \mathbf{M}_{1F})+O(\Delta t^2).
\end{eqnarray}
\end{subequations}

Let
\begin{equation}\label{eq:M-NCDE-4}
M_{0F}=S, A_{10}=-B_{00},
\end{equation}
then we can obtain $-\frac{B_{00}}{A_{10}}M_{0F}=S$, and Eq. (\ref{eq:ME-NCDE-2}) becomes
\begin{subequations}\label{eq:ME-NCDE-3}
\begin{equation}
\partial_t M_0+\nabla\cdot \mathbf{M}_1 = S+O(\Delta t),
\end{equation}
\begin{eqnarray}
&\partial_t M_0+\nabla\cdot \mathbf{M}_1+ \frac{\Delta t}{A_{10}}\big[A_{20}\partial_t^2 M_0+A_{21}\partial_t \nabla\cdot \mathbf{M}_1+A_{22}\nabla\nabla : \mathbf{M}_2\big]\nonumber\\
&=S-\frac{\Delta t}{A_{10}}(B_{10}\partial_t S+B_{11}\nabla\cdot \mathbf{M}_{1F})+O(\Delta t^2).
\end{eqnarray}
\end{subequations}
From Eq. (\ref{eq:ME-NCDE-3}a), we have
\begin{equation}\label{eq:ME-NCDE-3-1}
\partial_t^2 M_0=-\partial_t\nabla\cdot \mathbf{M}_1 + \partial_t S+O(\Delta t),
\end{equation}
then inserting Eq. (\ref{eq:ME-NCDE-3-1}) into Eq. (\ref{eq:ME-NCDE-3}b) yields
\begin{eqnarray}\label{eq:ME-NCDE-3-2}
&\partial_t M_0+\nabla\cdot \mathbf{M}_1+ \frac{\Delta t}{A_{10}}(B_{10}+A_{20})\partial_t S\nonumber\\
&=S-\frac{\Delta t}{A_{10}}\big[(A_{21}-A_{20})\partial_t \nabla\cdot \mathbf{M}_1+B_{11}\nabla\cdot \mathbf{M}_{1F}+A_{22}\nabla\nabla : \mathbf{M}_2\big]+O(\Delta t^2).
\end{eqnarray}

If we set
\begin{equation}\label{eq:M-NCDE-5}
(B_{10}+A_{20})\partial_t S=0,
\end{equation}
the moment equation at the leading order of $O(\Delta t^2)$ can be derived from Eq. (\ref{eq:ME-NCDE-3-2}),
\begin{eqnarray}\label{eq:ME-NCDE-3-2-1}
\partial_t M_0+\nabla\cdot \mathbf{M}_1=S-\frac{\Delta t}{A_{10}}\big[(A_{21}-A_{20})\partial_t \nabla\cdot \mathbf{M}_1+B_{11}\nabla\cdot \mathbf{M}_{1F}+A_{22}\nabla\nabla : \mathbf{M}_2\big]+O(\Delta t^2).
\end{eqnarray}

For a specific  NCDE, we can correctly recover the NCDE from Eq. (\ref{eq:ME-NCDE-3-2-1}) through choosing the moments and parameters. For example, taking
\begin{equation}\label{eq:M-NCDE-6}
M_0=\phi, ~\mathbf{M}_1=\mathbf{B},~\mathbf{M}_2=\beta c_s^2 \mathbf{D}+\mathbf{C}, ~\mathbf{M}_{1F}=-\frac{1}{B_{11}}\big[(A_{21}-A_{20})\partial_t \mathbf{B}+A_{22}\nabla\cdot \mathbf{\mathbf{C}}\big], ~\alpha=-\frac{A_{22}}{A_{10}}\Delta t \beta c_s^2 
\end{equation}
in Eq. (\ref{eq:ME-NCDE-3-2-1}) leads to the general NCDE \cite{Chai2020,ChaiYuanShi2023} with a constant diffusion coefficient $\alpha$
\begin{eqnarray}\label{eq:NCDE0}
\partial_t \phi+\nabla\cdot \mathbf{B}=S+\alpha\nabla\nabla : \mathbf{D}+O(\Delta t^2)=S+\nabla\cdot[\alpha\nabla\cdot \mathbf{D}]+O(\Delta t^2).
\end{eqnarray}

(2) Poisson-type equation or steady CDE \cite{LiChaiShi2014}

Let $M_0=0$ or consider a steady case, Eq. (\ref{eq:ME-NCDE-1}b) becomes
\begin{eqnarray}\label{eq:ME-NCDE-P}
A_{11}\nabla\cdot \mathbf{M}_1 +\Delta tA_{22}\nabla\nabla : \mathbf{M}_2+B_{00}M_{0F}+\Delta tB_{11}\nabla\cdot \mathbf{M}_{1F}=O(\Delta t^2).
\end{eqnarray}

Taking $\mathbf{M}_1=\mathbf{B}, ~\mathbf{M}_2=\beta c_s^2 \mathbf{D}, ~M_{0F}=S, ~\mathbf{M}_{1F}=\mathbf{0}, ~B_{00}=-A_{11}$, and $ \alpha=-\frac{A_{22}}{A_{11}}\Delta t \beta c_s^2$ in Eq. (\ref{eq:ME-NCDE-P}), the steady CDE can be recovered
\begin{eqnarray}\label{eq:ME-NCDE-P-1}
\nabla\cdot \mathbf{B} =\alpha\nabla\nabla : \mathbf{D}+S=O(\Delta t^2).
\end{eqnarray}

(3) Wave equation without source term.

Let $A_{20}=-1, ~A_{lm}=0$ for $(l,m)\neq (2,0) $ and $(2,2)$, and $B_{lm}=0$, for all $l,m$ in Eq. (\ref{eq:ME-NCDE-1}c) a wave equation can be recovered
\begin{eqnarray}\label{eq:ME-NCDE-W}
\partial_t^2 M_0 =A_{22}\nabla\nabla : \mathbf{M}_2+O(\Delta t^2).
\end{eqnarray}

Taking $M_0=\phi$ and $\mathbf{M}_2=\beta c_s^2 \mathbf{D}$ in Eq. (\ref{eq:ME-NCDE-W}), we obtain the following wave equation,
\begin{equation}\label{eq:ME-NCDE-W-1}
\partial_t^2 \phi =\alpha \nabla\nabla : \mathbf{D}+O(\Delta t^2)
\end{equation}
with $\alpha=A_{22}\beta c_s^2$.

\textbf{Remark 1}: (1) When $\mathbf{M}_1=\mathbf{B}=\mathbf{0}$, Eq. (\ref{eq:NCDE0}) becomes a diffusion equation, and we can choose $\mathbf{C}=0$, thus $ \mathbf{M}_{1F}=0$. In this case, there is no restriction on $A_{11}$, and only $A_{10}\neq 0$ is needed.

(2) If the source term $S=0$ or $S$ is a nonzero constant, then Eq. (\ref{eq:M-NCDE-5}) is satisfied. Otherwise, $B_{10}+A_{20}=0$ is needed.

(3) When $\mathbf{D}=\mathbf{0}$, Eq. (\ref{eq:NCDE0}) becomes a CE, which means that the CE can also be recovered directly from the MesoFD method. Another way to derive CE from the MesoFD method is to modify the diffusion coefficient $\alpha$ as a small quantity of order $\Delta t$ or $\Delta x$, \textcolor{blue}{such as Lax-Wendroff scheme.}

(4) If $\mathbf{M}_1=0$ and $\mathbf{M}_3=0$ in Eq. (\ref{eq:ME-NCDE-1}c), then to obtain Eq. (\ref{eq:ME-NCDE-W}), $A_{11}, A_{21}$ and $A_{33}$ can be non-zero.

(5) If $f_j(\mathbf{x},t)$ in Eq. (\ref{eq:2-1}) is replaced by $f_j(\mathbf{x},t+\Delta t)$, then $A_{l0}=\frac{1}{l!}\sum_{k,q}a_{kq}\cdot q^l$ ($ l\geq 1$) in Eqs. (\ref{eq:4-2}) and (\ref{eq:4-2-0}) for the moment equations (\ref{eq:ME-NCDE-0}) and (\ref{eq:ME-NCDE-1}) should be replaced by $A_{l0}=\frac{1}{l!}\big[\sum_{k,q}a_{kq}\cdot q^l-1\big]$  ($ l\geq 1$).

\subsection{The MesoFD and MMFD schemes with a three-level explicit structure for PDE}

As we know, the LB method for NCDE with a source term is commonly a three-level explicit scheme without correction term, or a two-level implicit scheme with the correction term. We now consider a three-level explicit scheme derived from Eq. (\ref{eq:2-1}),
\begin{equation}\label{eq:3layer}
 f_j(\mathbf{x},t+\Delta t)=\sum_{q=-1}^{0}\sum_{k=-m_q}^{m_q}\big[a_{kq}f_j^{eq[k,q]}+\Delta t b_{kq}F_j^{[k,q]}\big],
\end{equation}
and the macroscopic form (MMFDS) can be given as
\begin{equation}\label{eq:3layerMMFD}
 \phi(\mathbf{x},t+\Delta t)=\sum_j \sum_{q=-1}^{0}\sum_{k=-m_q}^{m_q}\big[a_{kq}f_j^{eq[k,q]}+\Delta t b_{kq}F_j^{[k,q]}\big].
\end{equation}

For simplicity, if we take $ k \in \{-1,0,1\}$, and let $a_{k0}=a_k, a_{k,-1}=\hat{a}_k, b_{k0}=b_k, b_{k,-1}=\hat{b}_k$, then Eqs. (\ref{eq:3layer}) and (\ref{eq:3layerMMFD}) become
\begin{eqnarray}\label{eq:3layer-1}
 f_j(\mathbf{x},t+\Delta t)&= & a_1 f_j^{eq}(\mathbf{x}+\mathbf{c}_j \Delta t,t)+ a_0 f_j^{eq}(\mathbf{x},t)+ a_{-1}f_j^{eq}(\mathbf{x}-\mathbf{c}_j \Delta t,t)\nonumber\\
  & &  + \hat{a}_1 f_j^{eq}(\mathbf{x}+\mathbf{c}_j \Delta t,t-\Delta t)+ \hat{a}_0 f_j^{eq}(\mathbf{x},t-\Delta t) +\hat{a}_{-1}f_j^{eq}(\mathbf{x}-\mathbf{c}_j \Delta t,t-\Delta t)\nonumber\\
   & &  + \Delta t  \big[b_1 F_j(\mathbf{x}+\mathbf{c}_j \Delta t,t)+ b_0 F_j(\mathbf{x},t)+ b_{-1}F_j(\mathbf{x}-\mathbf{c}_j \Delta t,t)\big]\nonumber\\
 & &    + \Delta t  \big[\hat{b}_1 F_j(\mathbf{x}+\mathbf{c}_j \Delta t,t-\Delta t)+ \hat{b}_0 F_j(\mathbf{x},t-\Delta t) +\hat{b}_{-1}F_j(\mathbf{x}-\mathbf{c}_j \Delta t,t-\Delta t)\big],
\end{eqnarray}
and
\begin{eqnarray}\label{eq:3layerMMFD-1}
 \phi(\mathbf{x},t+\Delta t)&=&  \sum_j\big[a_1 f_j^{eq}(\mathbf{x}+\mathbf{c}_j \Delta t,t)+ a_0 f_j^{eq}(\mathbf{x},t)+ a_{-1}f_j^{eq}(\mathbf{x}-\mathbf{c}_j \Delta t,t)\big]\nonumber\\
  & &   + \sum_j\big[\hat{a}_1 f_j^{eq}(\mathbf{x}+\mathbf{c}_j \Delta t,t-\Delta t)+ \hat{a}_0 f_j^{eq}(\mathbf{x},t-\Delta t) +\hat{a}_{-1}f_j^{eq}(\mathbf{x}-\mathbf{c}_j \Delta t,t-\Delta t)\big]\nonumber\\
  & &   + \Delta t  \sum_j\big[b_1 F_j(\mathbf{x}+\mathbf{c}_j \Delta t,t)+ b_0 F_j(\mathbf{x},t)+ b_{-1}F_j(\mathbf{x}-\mathbf{c}_j \Delta t,t)\big]\nonumber\\
   & &  + \Delta t  \sum_j\big[\hat{b}_1 F_j(\mathbf{x}+\mathbf{c}_j \Delta t,t-\Delta t)+ \hat{b}_0 F_j(\mathbf{x},t-\Delta t) +\hat{b}_{-1}F_j(\mathbf{x}-\mathbf{c}_j \Delta t,t-\Delta t)\big].
\end{eqnarray}

Using Eq. (\ref{eq:2-2-1}), we can rewrite Eq. (\ref{eq:3layerMMFD-1}) into another form,
\begin{eqnarray}\label{eq:3layerMMFD-2}
    \phi(\mathbf{x},t+\Delta t)=A_0 \phi(\mathbf{x},t)+\hat{A}_0\phi(\mathbf{x},t-\Delta t) +\Delta t\big[ B_0 S(\mathbf{x},t)+\hat{B}_0 S(\mathbf{x},t-\Delta t)\big]\nonumber\\
    + \frac{1}{2}\sum_j\big[a_1^{+}\Delta_{j}^{(2)}f_j^{eq,+}+\Delta t b_1^{+}\Delta_{j}^{(2)}F_j^{+}\big](\mathbf{x},t)
    + \frac{1}{2}\sum_j\big[\hat{a}_1^{+}\Delta_{j}^{(2)}f_j^{eq,+}+\Delta t \hat{b}_1^{+}\Delta_{j}^{(2)}F_j^{+}\big](\mathbf{x},t-\Delta t)\nonumber\\
    + \frac{1}{2}\sum_j\big[a_1^{-}\Delta_{j}^{(1)}f_j^{eq,-}+\Delta t b_1^{-}\Delta_{j}^{(1)}F_j^{-}\big](\mathbf{x},t)
    + \frac{1}{2}\sum_j\big[\hat{a}_1^{-}\Delta_{j}^{(1)}f_j^{eq,-}+\Delta t \hat{b}_1^{-}\Delta_{j}^{(1)}F_j^{-}\big](\mathbf{x},t-\Delta t),
\end{eqnarray}
where
\begin{eqnarray}\label{eq:2-7-1}
&A_0=a_0+a_1+a_{-1}, ~\hat{A}_0=\hat{a}_0+\hat{a}_1+\hat{a}_{-1}, ~B_0=b_0+b_1+b_{-1}, ~\hat{B}_0=\hat{b}_0+\hat{b}_1+\hat{b}_{-1},\nonumber\\
&a_1^{+}=a_1+a_{-1},a_1^{-}=a_1-a_{-1},b_1^{+}=b_1+b_{-1},b_1^{-}=b_1-b_{-1},\nonumber\\
&\hat{a}_1^{+}=\hat{a}_1+\hat{a}_{-1},\hat{a}_1^{-}=\hat{a}_1-\hat{a}_{-1},\hat{b}_1^{+}=\hat{b}_1+\hat{b}_{-1},\hat{b}_1^{-}=\hat{b}_1-\hat{b}_{-1}.
\end{eqnarray}

According to \textbf{Remark 1}, for Eqs. (\ref{eq:3layer-1}), (\ref{eq:3layerMMFD-1}) and (\ref{eq:3layerMMFD-2}), Eq. (\ref{eq:4-2-0}) has the following form,
\begin{subequations}\label{eq:4-2-00}
\begin{eqnarray}
&A_{00}=a_0 +a_1+ a_{-1}+\hat{a}_0 +\hat{a}_1+ \hat{a}_{-1}=1, 
~A_{10}=-(\hat{a}_0 +\hat{a}_1+ \hat{a}_{-1})-1, ~A_{11}=a_1-a_{-1}+\hat{a}_1-\hat{a}_{-1}, \nonumber\\
&A_{20}=\frac{1}{2}(\hat{a}_0 +\hat{a}_1+ \hat{a}_{-1}-1), ~A_{21}=-(\hat{a}_1-\hat{a}_{-1}),  ~A_{22}=\frac{1}{2}(a_1+a_{-1}+\hat{a}_1+ \hat{a}_{-1}),
\end{eqnarray}
\begin{eqnarray}
B_{00}=b_0 +b_1+ b_{-1}+\hat{b}_0 +\hat{b}_1+ \hat{b}_{-1}, ~
B_{10}=-(\hat{b}_0 +\hat{b}_1+ \hat{b}_{-1}), B_{11}=b_1-b_{-1}+\hat{b}_1-\hat{b}_{-1},
\end{eqnarray}
\end{subequations}
from which we can obtain the relation $A_{20}=-A_{10}/2-1$.

We note that Eq. (\ref{eq:3layerMMFD-1}) or (\ref{eq:3layerMMFD-2}) is a general form of MMFD method with a three-level structure, and is also a unified framework of FD method for convection-diffusion-like equations in $d$ dimensional space. Some existed FD schemes derived from the LBM can be taken as its special cases through selecting the parameters $a_k, \hat{a}_k, b_k$ and $\hat{b}_k$.

(1) SRT LB model with $\tau=1$: $a_0=a_1=0, a_{-1}=1; \hat{a}_0=\hat{a}_1=\hat{a}_{-1}=0$; $b_0=b_1=\hat{b}_0=\hat{b}_1=0,b_{-1}=3/2,\hat{b}_{-1}=-1/2$.

(2) TRT LB model with $s^{-}+s^{+}=2$ \cite{dHumieres2009,Ginzburg2012,Ginzburg2010}: $a_0=0,a_1=1-s^{-},a_{-1}=1; \hat{a}_0=-(1-s^{-}),\hat{a}_1=\hat{a}_{-1}=0$; $b_0=\hat{b}_0=-(1-s^{-})/2,b_1=\hat{b}_1=0,b_{-1}=3/2,\hat{b}_{-1}=-1/2$.

(3) One-stage simplified LB model without source term \cite{Gutierrez2021,RosisLiu2021}: $a_0=2(1-\tau),a_1=\tau-1, a_{-1}=\tau; \hat{a}_0=\hat{a}_1=\hat{a}_{-1}=0$.

(4) A new FD scheme for the wave equation without source term (Eq. (\ref{eq:ME-NCDE-W-1})): $a_0=2(1-\gamma), a_1=a_{-1}=\gamma; \hat{a}_0=-1, \hat{a}_1=\hat{a}_{-1}=0$, where $\gamma= \alpha/(\beta c_s^2)$.

In general, to derive a FD scheme, one can determine the parameters with Eq. (\ref{eq:4-2-00}). We now show that Eq. (\ref{eq:4-2-00}a) is a solvable linear system. In fact, let $X=(a_1, a_{-1}, \hat{a}_0, \hat{a}_{-1}, \hat{a}_{-1})^{T}$, $b=(1-a_0, 1+A_{10}, A_{11},A_{21},A_{22})^{T}$, then Eq.(\ref{eq:4-2-00}a) can be transformed into a solvable linear system $AX=b$ with the invertible coefficient matrix $A$,
\begin{equation}\label{eq:AX_b}
A=\left(
\begin{array}{ccccc}
    1 & 1  & 1 &  1 &  1  \\
    0 & 0  & -1 & -1 & -1  \\
    1 & -1  & 0 &  1 & -1  \\
    0 & 0  & 0 &  -1 &  1  \\
    1/2 & 1/2  & 0 &  1/2 &  1/2  \\
\end{array}
\right),
\end{equation}
and the inverse $A^{-1}$ can be obtained
\begin{equation}\label{eq:invA}
A^{-1}=\frac{1}{2}\left(
\begin{array}{ccccc}
    1 & 1  & 1 &  1 &  0  \\
    1 & 1  & -1 & -1 & 0  \\
    2 & 0  & 0 &  0 & -4  \\
    -1 & -1  & 0 &  -1 &  2  \\
    -1 & -1  & 0 &   1 &  2  \\
\end{array}
\right),
\end{equation}
therefore, we have
\begin{equation}\label{eq:X}
\left(
\begin{array}{c}
a_1\\
a_{-1}\\
\hat{a}_0\\
\hat{a}_1\\
\hat{a}_{-1}\\
\end{array}
\right)=
X=A^{-1}b=\frac{1}{2}\left(
\begin{array}{c}
-a_0+A_{10}+A_{11}+A_{21}+2\\
-a_0+A_{10}-A_{11}-A_{21}+2\\
-2a_0-4A_{22}+2\\
a_0-A_{10}-A_{21}+2A_{22}-2\\
a_0-A_{10}+A_{21}+2A_{22}-2\\
\end{array}
\right),
\end{equation}
with $A_{20}=-A_{10}/2-1$, while $\{b_k,\hat{b}_k, k=-1,0,1\}$ are of multiple selection.

\subsection{The stability of mesoscopic numerical method based macroscopic finite difference scheme for linear CDE and WE}

In this part, we consider the stability of MMFD scheme (\ref{eq:2-2}) for the linear CDE (LCDE) and WE (LWE) without the source term, respectively. The LCDE and LWE have the following forms,
\begin{eqnarray}\label{eq:LCDE}
\partial_t \phi+\nabla\cdot \mathbf{u}\phi=\nabla \cdot\big[\alpha \nabla \cdot\mathbf{D}_0\phi\big],
\end{eqnarray}
\begin{eqnarray}\label{eq:LWE}
\partial_t^2 \phi=\nabla \cdot\big[\alpha \nabla \cdot\mathbf{D}_0\phi\big],
\end{eqnarray}
where $\mathbf{u}$ is a constant vector, and $\mathbf{D}_0$ is a constant matrix (tensor). When $\mathbf{D}_0=\mathbf{0}$, the LCDE will become a CE. The MesoFD scheme for Eq. (\ref{eq:LCDE}) and Eq. (\ref{eq:LWE}) can be simplified as
\begin{equation}\label{eq:LMMFD}
f_j(\mathbf{x},t)=\sum_{k,q}a_{kq}f_j^{eq[k,q]}.
\end{equation}

Let $\hat{f}(t)$ be the Fourier transformation of a function $f(\mathbf{x},t)$, then we have
\begin{equation}\label{eq:TforMMFD}
\hat{f}_j(t)=\sum_{k,q}a_{kq}E_j^k\hat{f}_j^{eq}(t+q\Delta t), E_j=\exp[\textbf{i} \xi_j ], \xi_j=\Delta t \mathbf{\mathbf{c}}_j\cdot\xi.
\end{equation}
For Eqs. (\ref{eq:LCDE}) and (\ref{eq:LWE}), $f^{eq}_j$ can be taken the same form \cite{ChaiYuanShi2023},
\begin{equation}\label{eq:TforMMFD-0}
f_j^{eq}=\phi \bar{u}_j, \bar{u}_j=\omega_j\left(1+\frac{\mathbf{c}_j\cdot \mathbf{u}}{c_s^2}+\frac{\mathbf{C}_{\alpha\alpha}\mathbf{Q}_{j\alpha\alpha}}{c_s^2(c_{\alpha}^2-c_s^2)}+\frac{\mathbf{C}_{\alpha\bar{\alpha}}\mathbf{Q}_{j\alpha\bar{\alpha}}}{2c_s^4}\right),
\end{equation}
where $\omega_j$ is the weight coefficient (see Appendix A for details),  $\mathbf{C}=c_s^2(\mathbf{\beta D}_0-\mathbf{I})+\bar{\lambda}\mathbf{uu}, ~\mathbf{Q}_j =\mathbf{c}_j\mathbf{c}_j-c_s^2 \mathbf{I}, ~c_s^2=d_{0\alpha} c_{\alpha}^2, ~c_{\alpha}=\Delta x_{\alpha}/\Delta t$, and $\bar{\lambda} =0$ or 1. For the WE (\ref{eq:LWE}), we can set $\mathbf{u}=0$ in Eq. (\ref{eq:TforMMFD-0}). Then summing Eq. (\ref{eq:TforMMFD}), we obtain
\begin{equation}\label{eq:TforMMFD-1}
\hat{\phi}(t)=\sum_{k,q}a_{kq}\sum_j E_j^k \bar{u}_j\hat{\phi}(t+q\Delta t).
\end{equation}
Let
\begin{equation}\label{eq:TforMMFD-2}
\bar{u}_j=\bar{u}_j^{+}+\bar{u}_j^{-},~\bar{u}_j^{+}=\omega_j\left(1+\frac{\mathbf{C}_{\alpha\alpha}\mathbf{Q}_{j\alpha\alpha}}{c_s^2(c_{\alpha}^2-c_s^2)}+\frac{\mathbf{C}_{\alpha\bar{\alpha}}\mathbf{Q}_{j\alpha\bar{\alpha}}}{2c_s^4}\right),~\bar{u}_j^{-}=\omega_j\frac{\mathbf{c}_j\cdot \mathbf{u}}{c_s^2},
~A_q=\sum_k a_{kq}\sum_j E_j^k \bar{u}_j,
\end{equation}
we can obtain
\begin{subequations}\label{eq:TforMMFD-3}
\begin{eqnarray}
\hat{\phi}(t)=\sum_{q=-n}^0 A_q\hat{\phi}(t+q\Delta t),
\end{eqnarray}
\begin{eqnarray}
A_q&=&\sum_k a_{kq} \sum_{j\geq 0}\big[\bar{u}_j^{+}\cos(k\xi_j)+\textbf{i}\bar{u}_j^{-}\sin(k\xi_j)\big]\nonumber\\
&=&\sum_k a_{kq}+ \sum_{k\geq 1} \sum_{j\geq 0}\big[a_{kq}^{+}\bar{u}_j^{+}(\cos(k\xi_j)-1)+\textbf{i}a_{kq}^{-}\bar{u}_j^{-}\sin(k\xi_j)\big],~ q\leq 0,
\end{eqnarray}
\end{subequations}
where $\sum_j \bar{u}_j^{+}=1$ has been used, $a_{kq}^{+} and a_{kq}^{-}$ are defined by Eq. (\ref{eq:2-1-FD1-0}).
Actually, Eq. (\ref{eq:TforMMFD-3}a) is equivalent to
\begin{equation}\label{eq:TforMMFD-4}
\hat{\phi}(t)=\sum_{q=-n}^{-1} \frac{A_q}{1-A_0}\hat{\phi}(t+q\Delta t).
\end{equation}

The characteristic polynomial of Eq. (\ref{eq:TforMMFD-4}) is of degree $n$, and has the following form
\begin{equation}\label{eq:TforMMFD-4-0}
p(\lambda)=\lambda^n-\sum_{k=0}^{n-1} p_{k-n} \lambda^k, ~p_q=\frac{A_q}{1-A_0},~q=-n,~..., ~-1,
\end{equation}
the MMFD scheme is stable if the modulus of any root $\lambda$ of this polynomial is no greater than 1. Usually, when $n$ is large, one cannot easily obtain the stability condition. For this reason, here we consider some special cases.

For a two-level scheme with $q \in \{0,-1\}$, Eq. (\ref{eq:TforMMFD-4}) can be simplified as
\begin{equation}\label{eq:TforMMFD-5}
\hat{\phi}(t)= \frac{A_{-1}}{1-A_0}\hat{\phi}(t-\Delta t),
\end{equation}
and the stability condition is given by
\begin{equation}\label{eq:TforMMFD-6}
|\lambda|=|\frac{A_{-1}}{1-A_0}|\leq 1.
\end{equation}

If $A_0=0$, then a two-level explicit scheme is obtained with
\begin{eqnarray}\label{eq:TforMMFD-7}
\lambda=A_{-1}=\sum_{k=-m}^{m} a_k \sum_{j\geq 0}\big[\bar{u}_j^{+}\cos(k\xi_j)+\textbf{i}\bar{u}_j^{-}\sin(k\xi_j)\big]\nonumber\\
=1+\sum_{j\geq 1} \sum_{k=1}^{m}\big[(a_k+a_{-k})\bar{u}_j^{+}(\cos(k\xi_j)-1)+\textbf{i}(a_k-a_{-k})\bar{u}_j^{-}\sin(k\xi_j)\big],
\end{eqnarray}
where $\sum_{k} a_{k,-1}=1$ is used, $a_k=a_{k,-1}, a_{-k}=a_{-k,-1}, k=1,...,m$, and $m=m_{-1}$.

Note that Eq. (\ref{eq:TforMMFD-7}) has the same form as Eq. (77) in Ref. \cite{Hindmarsh1984}. Based on Eq. (\ref{eq:TforMMFD-7}), the stability theorem in Ref. \cite{Hindmarsh1984} can be extended as follows.

\textbf{Theorem 1} The two-level explicit scheme for LCDE without the source term is stable ($|\lambda|\leq 1$) if and only if the following conditions holds,
\begin{subequations}\label{eq:Tcond}
\begin{eqnarray}
(a_k+a_{-k})\bar{u}_j^{+}>0, \forall k\geq 1, j\geq 1,
\end{eqnarray}
\begin{eqnarray}
\sum_{j\geq 1}\sum_{k=1}^{m}(a_k+a_{-k})\bar{u}_j^{+}\leq 1,
\end{eqnarray}
\begin{eqnarray}
\sum_{j\geq 1}\sum_{k=1}^{m}\frac{[(a_k-a_{-k})\bar{u}_j^{-}]^2}{(a_k+a_{-k})\bar{u}_j^{+}}\leq 1.
\end{eqnarray}
\end{subequations}

{\textbf{Remark 2} It is found from Eq. (\ref{eq:Tcond}c) that $ [(a_k-a_{-k})\bar{u}_j^{-}]^2 \leq (a_k+a_{-k})\bar{u}_j^{+}$, even if $(a_k+a_{-k})\bar{u}_j^{+}=0$. If any $(a_k+a_{-k})\bar{u}_j^{+}=0$, then $(a_k-a_{-k})\bar{u}_j^{-}=0$, and thus Eq. (\ref{eq:Tcond}c) follows by summing only over those $j$ and $k$ for which $(a_k+a_{-k})\bar{u}_j^{+}>0$. 

Since
\begin{equation}\label{eq:Tcond-1}
A_{00}=a_0+\sum_{k=1}^{m}(a_k+a_{-k})=1,~\sum_j\bar{u}_j^{+}= 1, ~\sum_{j\geq 0}\big[a_0+\sum_{k=1}^{m}(a_k+a_{-k})\big]\bar{u}_j^{+}= 1,
\end{equation}
we have
\begin{eqnarray}\label{eq:Tcond-2}
\sum_{j\geq 1}\sum_{k=1}^{m}(a_k+a_{-k})\bar{u}_j^{+}= 1-(a_0+\bar{u}_0^{+}(1-a_0))\leq 1 \Leftrightarrow a_0+\bar{u}_0^{+}(1-a_0)\geq 0\nonumber\\
\Leftrightarrow (1-a_0)(1- \bar{u}_0^{+})\leq 1 \Leftrightarrow (1-\bar{u}_0^{+})\sum_{k=1}^{m}(a_k+a_{-k})\leq 1.
\end{eqnarray}
Eq. (\ref{eq:Tcond}b) is equivalent to
\begin{equation}\label{eq:Tcond-3}
a_0+\bar{u}_0^{+}(1-a_0)\geq 0.
\end{equation}

For the two level $\theta$-scheme (Eq. (\ref{eq:GCN})) without the source term, one can obtain
\begin{equation}\label{eq:Tcond1-1}
\lambda=\frac{1+\theta (A+\textbf{i}B)}{1-(1-\theta)(A+\textbf{i}B)},
\end{equation}
where
\begin{equation}\label{eq:Tcond1-2}
A=\sum_j\sum_{k=1}^{m_0}(a_k+a_{-k})\bar{u}_j^{+}(\cos (k\xi _j)-1),B=\sum_j\sum_{k=1}^{m_0}(a_k-a_{-k})\bar{u}_j^{-}\sin (k\xi _j),
\end{equation}
and $\sum_k a_k=1$ is used. Then we obtain the following equation from Eqs. (\ref{eq:Tcond1-1}) and (\ref{eq:Tcond1-2})
\begin{equation}\label{eq:Tcond1-3}
|\lambda|^2-1=\frac{2A+(2\theta -1)(A^2+B^2)}{|1-(1-\theta)(A+iB)|^2}.
\end{equation}
Therefore, if Eq. (\ref{eq:Tcond}a) holds, we can prove that $|\lambda|\leq 1$ for $0\leq \theta \leq 1/2 $ due to $A\leq 0$, while the conditional stability is obtained for $\theta > 1/2$. Based above analysis, we have the stability condition of the two level $\theta$-scheme (Eq. (\ref{eq:GCN})) without the source term.

\textbf{Theorem 2} The two-level $\theta$-scheme (Eq. (\ref{eq:GCN})) for LCDE without the source term is (1) unconditional stable  ($|\lambda|\leq 1$) for $0\leq \theta \leq 1/2 $, (2) conditional stable for $\theta > 1/2$ if Eq. (\ref{eq:Tcond}a) holds. 

Further, taking $n=2$ in Eq. (\ref{eq:TforMMFD-4-0}) we can obtain the following characteristic polynomial of a general three-level scheme,
\begin{equation}\label{eq:3-level}
p(\lambda)=\lambda^2-p_1 \lambda-p_0,p_0=\frac{A_{-2}}{1-A_0},p_1=\frac{A_{-1}}{1-A_0},
\end{equation}
where $A_0, A_{-1} $, and $A_{-2}$ are defined by Eq. (\ref{eq:TforMMFD-3}b). Then one can derive the following result from Schur-Cohn theorem \cite{Miller1971}.

{\textbf{Theorem 3} The three-level MMFD scheme is stable if and only if the following conditions holds,
\begin{equation}\label{eq:3-level-1}
|p_1+ \bar{p}_{1} p_0|\leq 1-|p_0|^2,
\end{equation}
where $\bar{p}_1$ is the conjugate of $p_1$.

Some corollaries of \textbf{Theorem 3} can be obtained for the specific three-level MMFD scheme. In the following, we take $A_0=0$ in Eq. (\ref{eq:3-level}), and two explicit schemes are considered.

For the TRT-LB model with $s^{+}+s^{-}=2$ as mentioned above, we have
$p_0=s^{-}-1,~p_1=(2-s^{-}) A-\textbf{i}s^{-} B$, where
\begin{equation}\label{eq:3-level-2}
A=\sum_j\bar{u}_j^{+}\cos (\xi _j),B=\sum_j\bar{u}_j^{-}\sin (\xi _j),
\end{equation}
then Eq. (\ref{eq:3-level-1}) becomes $|A-\textbf{i} B| \leq 1$, and we obtain

\textbf{Corollary 1} The TRT LB model with $s^{+}+s^{-}=2$ is stable if and only if 
\begin{equation}\label{eq:3-level-3}
A^2+B^2 \leq 1.
\end{equation}

A sufficient condition of Eq. (\ref{eq:3-level-3}) was given in Ref. \cite{Ginzburg2010}, which leads to the following corollary.

\textbf{Corollary 2} The TRT LB model with $s^{+}+s^{-}=2$ is stable if the following conditions hold, 
\begin{equation}\label{eq:3-level-4}
(1)~ \bar{u}_0^{+} \geq 0, \bar{u}_j^{+} >0, \forall j>0.~~ (2)~ \sum_{j\geq 1}\frac{(\bar{u}_j^{-})^2}{\bar{u}_j^{+}} \leq 1.
\end{equation}

{\textbf{Remark 3} Similarly to \textbf{Remark 2 }, if $\bar{u}_j^{+}=0$, then $\bar{u}_j^{-}=0$, and thus Eq. (\ref{eq:3-level-4}) follows by summing only over those $j$ for which $\bar{u}_j^{+}>0$.

For the three-level MMFD for WE mentioned above, we have
$p_0=-1,~p_1=2- \bar{A}$, where
\begin{equation}\label{eq:3-level-5}
\bar{A}=2\gamma \sum_j\bar{u}_j^{+}(1-\cos (k\xi _j))=4\gamma \sum_j\bar{u}_j^{+}\sin^2 (k\xi _j/2),~\gamma = \frac{\alpha}{\beta c_s^2},
\end{equation}
then Eq. (\ref{eq:3-level-1}) is always true.

Now, to simplify the conditions in Eqs. (\ref{eq:Tcond}) and (\ref{eq:3-level-4}), we consider the linear EDF. Taking $\mathbf{D}_0=\mathbf{I}, \beta=1, \bar{\lambda}=0$ in Eq. (\ref{eq:TforMMFD-0}), the linear EDF $f_j^{eq}$ with $\bar{u}_j^{+}=\omega_j >0$ and $\bar{u}_j^{-}=\omega_j\frac{\mathbf{c}_j\cdot \mathbf{u}}{c_s^2}$ is obtained, and
\begin{equation}\label{eq:Tcond-4}
\sum_{j\geq 1}\frac{(\bar{u}_j^{-})^2}{\bar{u}_j^{+}}=\sum_{j\geq 1} \omega_j\frac{(\mathbf{c}_j\cdot \mathbf{\mathbf{u}})^2}{c_s^4}=\frac{1}{c_s^4}\sum_{j\geq 1} \omega_j\mathbf{c}_j\mathbf{c}_j:\mathbf{uu}=\frac{|\mathbf{u}|^2}{c_s^2},
\end{equation}
where  $\sum_j \omega_j \mathbf{c}_j\mathbf{c}_j=c_s^2 \mathbf{I}$ is used. From Eq. (\ref{eq:Tcond-4}) it follows that Eq. (\ref{eq:Tcond}) and Eq. (\ref{eq:Tcond-3}) become
\begin{subequations}\label{eq:Tcond-LEDF}
\begin{eqnarray}
(a_k+a_{-k})>0, \forall k\geq 1,
\end{eqnarray}
\begin{eqnarray}
(1-\omega_0)(1-a_0)=(1-\omega_0)\sum_{k=1}^{m}(a_k+a_{-k})\leq 1,
\end{eqnarray}
\begin{eqnarray}
\frac{|\mathbf{u}|^2}{c_s^2}\sum_{k=1}^{m}\frac{(a_k-a_{-k})^2\emph{}}{(a_k+a_{-k})}\leq 1,
\end{eqnarray}
\end{subequations}
and
\begin{equation}\label{eq:Tcond-LEDF-1}
a_0+\omega_0(1-a_0)\geq 0.
\end{equation}

The equivalent form of Eq. (\ref{eq:LMMFD}) can be further derived from Eq. (\ref{eq:2-2-1}) as
\begin{eqnarray}\label{eq:LMMFD-1}
    \phi(\mathbf{x},t)=\sum_{k,q}a_{kq} \phi(\mathbf{x},t+q\Delta t)
    + \frac{1}{2}\sum_q\sum_j\sum_{k=1}^{m_q}(a_{kq}+a_{-kq})\bar{u}_j^{+}\Delta_{j,k}^{(2)}\phi(\mathbf{x},t+q\Delta t)\nonumber\\
+  \frac{1}{2}\sum_q\sum_j\sum_{k=1}^{m_q}(a_{kq}-a_{-kq})\bar{u}_j^{-}\Delta_{j,k}^{(1)}\phi(\mathbf{x},t+q\Delta t).
\end{eqnarray}
For the two-level explicit scheme, Eq. (\ref{eq:LMMFD-1}) can be simplified by
\begin{eqnarray}\label{eq:LMMFD-2}
    \phi(\mathbf{x},t)=\phi(\mathbf{x},t-\Delta t) + \frac{1}{2}\sum_j\sum_{k=1}^{m}\big[(a_{k}+a_{-k})\bar{u}_j^{+}\Delta_{j,k}^{(2)}+ (a_{k}-a_{-k})\bar{u}_j^{-}\Delta_{j,k}^{(1)}\big]\phi(\mathbf{x},t-\Delta t),
\end{eqnarray}
with
\begin{subequations}\label{eq:LMMFD-2-0}
\begin{eqnarray}
A_{10}=-\sum_{k=-m}^m a_{k}=-A_{00}=-1, A_{11}=\sum_{k=1}^m (a_{k}-a_{-k})=A_{10}=-1, A_{22}=\frac{1}{2}\sum_{k=1}^m (a_k+a_{-k}),
\end{eqnarray}
\begin{eqnarray}
\alpha= -\frac{A_{22}}{A_{10}}  c_s^2 \Delta t=\frac{ c_s^2 \Delta t}{2} \sum_{k=1}^{m}(a_k+a_{-k})\leq \frac{ c_s^2 \Delta t}{2(1-\omega_0)},
\end{eqnarray}
\end{subequations}
where Eq. (\ref{eq:4-2-0}), Eq. (\ref{eq:M-NCDE-6}) and Eq. (\ref{eq:Tcond-LEDF}b) have been used. Note that Eq. (\ref{eq:LMMFD-2-0}b) is equivalent to Eq. (\ref{eq:Tcond-LEDF}b) or (\ref{eq:Tcond-LEDF-1}). When $m=1$, Eq. (\ref{eq:LMMFD-2-0}a) is simplified as
\begin{equation}\label{eq:LMMFD-2-1}
A_{00}=a_0+a_1+a_{-1}=-A_{10}=1, A_{10}=A_{11}=a_1-a_{-1}=-1, A_{22}=\frac{1}{2} (a_1+a_{-1})=\frac{\alpha}{ \Delta t c_s^2},
\end{equation}
and the stability condition in Eq. (\ref{eq:Tcond-LEDF}) becomes 
\begin{subequations}\label{eq:Tcond-LEDF-3}
\begin{eqnarray}
a_1+a_{-1}=1-a_0=\frac{2\alpha}{\Delta t c_s^2}>0,
\end{eqnarray}
\begin{eqnarray}
(1-\omega_0)(1-a_0)=(1-\omega_0)(a_1+a_{-1})=(1-\omega_0)\frac{2\alpha}{ \Delta t c_s^2}\leq 1,
\end{eqnarray}
\begin{eqnarray}
\frac{|\mathbf{u}|^2}{c_s^2}\leq a_1+a_{-1}=\frac{2\alpha}{ \Delta t c_s^2},
\end{eqnarray}
\end{subequations}
then one can obtain from Eq. (\ref{eq:LMMFD-2-1}) that
\begin{equation}\label{eq:LMMFD-2-2}
a_0=1-\frac{2\alpha}{ \Delta t c_s^2}, ~a_1= \frac{1}{2}(\frac{2\alpha}{ \Delta t c_s^2}-1), ~a_{-1}=\frac{1}{2}(\frac{2\alpha}{ \Delta t c_s^2}+1),
\end{equation}
and the stability condition in Eq. (\ref{eq:Tcond-LEDF-3}) is equivalent to
\begin{equation}\label{eq:LMMFD-2-3}
|\mathbf{u}|^2 \leq \frac{2\alpha}{ \Delta t }\leq \frac{c_s^2}{1-\omega_0}.
\end{equation}

It is clear that the difference of stability region given by Eq. (\ref{eq:LMMFD-2-3})  depends on two parameters $c_s^2$ and $\omega_0$ in rD$d$Q$q$ lattice model.

For a rD$d$Q$b$ ($b=2\bar{m}+1$) lattice model with the discrete velocity set $\mathbf{V}_b=\{\mathbf{0},\mathbf{c}_1,...,\mathbf{c}_{\bar{m}},-\mathbf{c}_1,...-\mathbf{c}_{\bar{m}}\}$, Eq. (\ref{eq:LMMFD-2}) becomes
\begin{eqnarray}\label{eq:LMMFD-3}
    \phi(\mathbf{x},t)=\phi(\mathbf{x},t-\Delta t) + \sum_{j=1}^{\bar{m}}\sum_{k=1}^{m}\big[(a_{k}+a_{-k})\bar{u}_j^{+}\Delta_{j,k}^{(2)}+ (a_{k}-a_{-k})\bar{u}_j^{-}\Delta_{j,k}^{(1)}\big]\phi(\mathbf{x},t-\Delta t).
\end{eqnarray}
When $m=1$, Eq. (\ref{eq:LMMFD-3}) becomes
\begin{eqnarray}\label{eq:LMMFD-3-1}
    \phi(\mathbf{x},t)&=&\phi(\mathbf{x},t-\Delta t) + \sum_{j=1}^{\bar{m}}\big[(a_{1}+a_{-1})\bar{u}_j^{+}\Delta_{j}^{(2)}+ (a_{1}-a_{-1})\bar{u}_j^{-}\Delta_{j}^{(1)}\big]\phi(\mathbf{x},t-\Delta t)\nonumber\\
   & =&\phi(\mathbf{x},t-\Delta t) + \sum_{j=1}^{\bar{m}}\left[\frac{2\alpha \bar{u}_j^{+}}{ \Delta t c_s^2}\Delta_{j}^{(2)} -\bar{u}_j^{-}\Delta_{j}^{(1)}\right]\phi(\mathbf{x},t-\Delta t).
\end{eqnarray}

It should be noted that Eq. (\ref{eq:LMMFD-3-1}) is a general two-level FD scheme for LCDE (\ref{eq:LCDE}) based on the rectangular rD$d$Q$b$ lattice model \cite{ChaiYuanShi2023} and linear EDF. The FD scheme in Ref. \cite{Hindmarsh1984} for Eq. (\ref{eq:LCDE}) is a special case of the rectangular rD$d$Q$b$ lattice model \cite{ChaiYuanShi2023} with $\bar{m}=d $ and $\bar{\lambda}=0$. Since $ d_{0j}=\frac{c_s^2}{c_{j}^2},c_j= \frac{\Delta x_j}{\Delta t},\omega_{j}=\frac{d_{0j}}{2},j=1,...,d,\omega_0=1-2(\omega_1+...+\omega_d)$, it follows from Eqs. (\ref{eq:LMMFD-2-0}b) and (\ref{eq:LMMFD-3-1}) that
\begin{subequations}\label{eq:LMMFD-4}
\begin{eqnarray}
    \phi(\mathbf{x},t)=\phi(\mathbf{x},t-\Delta t) + \sum_{j=1}^{d}\left[\frac{\alpha \Delta t}{\Delta x_j^2}(\phi_j-2\phi+\phi_{\bar{j}}) - \frac{u_j\Delta t}{2 \Delta x_j}(\phi_j-\phi_{\bar{j}})\right],
\end{eqnarray}
\begin{eqnarray}
    \sum_{j=1}^d\frac{2\alpha \Delta t}{\Delta x_j^2}\leq \sum_{j=1}^d \frac{d_{0j}}{1-\omega_0} = \sum_{j=1}^d\frac{d_{0j}}{d_{01}+...+d_{0d}}=1,
\end{eqnarray}
\end{subequations}
where $\phi=\phi(\mathbf{x},t-\Delta t),\phi_j=\phi(\mathbf{x+c}_j\Delta t,t-\Delta t),\phi_{\bar{j}}=\phi(\mathbf{x-c}_j\Delta t,t-\Delta t)$. If we take rD$d Qb$ lattice model with $\mathbf{V}_{2^d+1}=\{\mathbf{0},(\pm 1,\pm 1,\cdot\cdot\cdot,\pm 1) \}$ and $b=2^d+1$, then $\omega_j=\omega_1, j>1, \omega_0=1-2^d \omega_1, d_0=2^d \omega_1$, which gives a square lattice ($\Delta x_j=\Delta x, \forall j$), and Eq. (\ref{eq:LMMFD-2-0}b) becomes
\begin{equation}\label{eq:LMMFD-2-4}
\sum_{j=1}^d\frac{2\alpha\Delta t}{\Delta x_j^2}\leq \sum_{j=1}^d\frac{d_0 }{1-\omega_0}=d.
\end{equation}

One can find that the condition in Eq. (\ref{eq:LMMFD-4}b) is much stronger than that in Eq. (\ref{eq:LMMFD-2-4}) for a large $d$.

\subsection{The derivation of the moment equations for NSEs using DTE method}

Now using the same way as above, we can derive the moment equation for NSEs from the MesoFD method (\ref{eq:2-1}). The basic moments of $f_j$, $f_j^{eq}$, and $F_j$ are given as follows,
\begin{subequations}\label{eq:M-NSEs}
\begin{equation}
M_0=\sum_j f_j^{eq}=\sum_j f_j,\ \  \mathbf{M}_1=\sum_j \mathbf{c}_j f_j^{eq}=\sum_j \mathbf{c}_j f_j,
\end{equation}
\begin{equation}
\mathbf{M}_2=\sum_j \mathbf{c}_j \mathbf{c}_j f_j^{eq}, \mathbf{M}_3=\sum_j \mathbf{c}_j \mathbf{c}_j \mathbf{c}_j f_j^{eq}, \mathbf{M}_k=\sum_j \mathbf{c}_j\cdot\cdot\cdot \mathbf{c}_j f_j^{eq}, ~k>3
\end{equation}
\begin{equation}
M_{0F}=\sum_j F_j,\ \ \mathbf{M}_{1F}=\sum_j \mathbf{c}_j F_j,\ \ \mathbf{M}_{2F}=\sum_j \mathbf{c}_j \mathbf{c}_j F_j, \mathbf{M}_{kF}=\sum_j \mathbf{c}_j\cdot\cdot\cdot \mathbf{c}_j F_j, ~k>2
\end{equation}
\end{subequations}
where $\mathbf{M}_k$ and $\mathbf{M}_{kF}$ ($k\geq 0$) are the $k$-th moments of $f_j, f_j^{eq}$ and $F_j$, respectively.
Additionally, Eq. (\ref{eq:M-NSEs}a) gives the following moments of nonequilibrium distribution function,
\begin{equation}\label{eq:M-NSEs-1}
M_0^{ne}=\sum_j f_j^{ne}=0,\ \ \mathbf{M}_1^{ne}=\sum_j \mathbf{c}_j f_j^{ne}=0.
\end{equation}

Similarly, using DTE method we can derive the moment equations from Eq.(\ref{eq:2-1}) or Eq. (\ref{eq:2-3}). Before doing that, we first use the double-distribution based scheme  (\ref{eq:2-3}) with the form of moments:
\begin{subequations}\label{eq:MMFD-NSEs}
\begin{equation}
M_0(\mathbf{x},t)=\sum_j\sum_{k,q}\big[a_{kq}f_j^{eq[k,q]}+\Delta t b_{kq}F_j^{[k,q]}\big],  \ \
\end{equation}
\begin{equation}
\mathbf{M}_1(\mathbf{x},t)=\sum_j\sum_{k,q}\mathbf{c}_j\big[\bar{a}_{kq} f_j^{eq[k,q]}+\Delta t \bar{b}_{kq} F_j^{[k,q]}\big],  \ \
\end{equation}
\end{subequations}
in which two different set of parameters, $\{a_{kq},b_{kq}\}$ and $\{\bar{a}_{kq},\bar{b}_{kq}\}$, are used.

Using Taylor expansion to Eq.(\ref{eq:MMFD-NSEs}), we have
\begin{subequations}\label{eq:ME-NSEs}
\begin{equation}
M_0=
A_{00}M_0+\sum_{l=1}^{N}\Delta t^l\sum_{m=0}^{l}\big[A_{lm}(\partial_t^{l-m} \nabla^m \otimes\mathbf{M}_m)\big]+ \sum_{l=0}^{N-1}\Delta t^{l+1}\sum_{m=0}^{l}\big[B_{lm}(\partial_t^{l-m} \nabla^m \otimes\mathbf{M}_{mF})\big]+O(\Delta t^{N+1}),
\end{equation}
\begin{equation}
\mathbf{M}_1=
\bar{A}_{00}\mathbf{M}_1+\sum_{l=1}^{N}\Delta t^l\sum_{m=0}^{l}\big[\bar{A}_{lm}(\partial_t^{l-m} \nabla^m \otimes\mathbf{M}_{m+1})\big]+ \sum_{l=0}^{N-1}\Delta t^{l+1}\sum_{m=0}^{l}\big[\bar{B}_{lm}(\partial_t^{l-m} \nabla^m \otimes\mathbf{M}_{{m+1}F})\big]+O(\Delta t^{N+1}), N\geq 1,
\end{equation}
\end{subequations}
where
\begin{equation}\label{eq:ME-NSEs-0}
\bar{A}_{lm}=\frac{1}{l!}C_l^m  \sum_{k,q}\bar{a}_{kq}q^{l-m}k^m,  \bar{B}_{lm}=\frac{1}{l!}C_l^m  \sum_{k,q}\bar{b}_{kq}q^{l-m}k^m,
\end{equation}
and $A_{lm}$ is defined as in Eq. (\ref{eq:4-2}).

Taking $\bar{A}_{00}=A_{00}=1$ in Eq. (\ref{eq:ME-NSEs}), we obtain a general moment equation for NSEs,
\begin{subequations}\label{eq:ME-NSEs-1}
\begin{equation}
\sum_{l=1}^{N}\Delta t^{l-1}\sum_{m=0}^{l}\big[A_{lm}(\partial_t^{l-m} \nabla^m \otimes\mathbf{M}_m)\big]+ \sum_{l=0}^{N-1}\Delta t^{l}\sum_{m=0}^{l}\big[B_{lm}(\partial_t^{l-m} \nabla^m \otimes\mathbf{M}_{mF})\big]=O(\Delta t^{N}),
\end{equation}
\begin{equation}
\sum_{l=1}^{N}\Delta t^{l-1}\sum_{m=0}^{l}\big[\bar{A}_{lm}(\partial_t^{l-m} \nabla^m \otimes\mathbf{M}_{m+1})\big]+ \sum_{l=0}^{N-1}\Delta t^{l}\sum_{m=0}^{l}\big[\bar{B}_{lm}(\partial_t^{l-m} \nabla^m \otimes\mathbf{M}_{{m+1}F})\big]+O(\Delta t^{N}), N\geq 1.
\end{equation}
\end{subequations}

Note that Eq. (\ref{eq:ME-NSEs-1}a) is same as Eq. (\ref{eq:ME-NCDE-0}). In what follows, we only discuss Eq.(\ref{eq:ME-NSEs-1}b) , and  have the following equations for $N=1,2$,
\begin{subequations}\label{eq:ME-NSEs-2}
\begin{eqnarray}
\bar{A}_{10}\partial_t \mathbf{M}_1+\bar{A}_{11}\nabla\cdot \mathbf{M}_2 =-\bar{B}_{00}\mathbf{M}_{1F}+O(\Delta t),
\end{eqnarray}
\begin{eqnarray}
\bar{A}_{10}\partial_t \mathbf{M}_1+\bar{A}_{11}\nabla\cdot \mathbf{M}_2 + \Delta t\big[\bar{A}_{20}\partial_t^2 \mathbf{M}_1+\bar{A}_{21}\partial_t \nabla\cdot \mathbf{M}_2 +\bar{A}_{22}\nabla\nabla : \mathbf{M}_3\big]\nonumber\\
=-\bar{B}_{00}\mathbf{M}_{1F}-\Delta t(\bar{B}_{10}\partial_t \mathbf{M}_{1F}+\bar{B}_{11}\nabla\cdot \mathbf{M}_{2F} )+O(\Delta t^2).
\end{eqnarray}
\end{subequations}

Taking
\begin{subequations}\label{eq:M-NSEs-2}
\begin{equation}
A_{10}=A_{11}=-B_{00}, -\frac{B_{00}}{A_{10}}M_{0F}=S,
\end{equation}
\begin{equation}
\bar{A}_{10}=\bar{A}_{11}=-\bar{B}_{00}, -\frac{\bar{B}_{00}}{\bar{A}_{10}}\mathbf{M}_{1F}=\mathbf{F},
\end{equation}
\end{subequations}
Eq. (\ref{eq:ME-NSEs-2}) becomes
\begin{subequations}\label{eq:ME-NSEs-3}
\begin{eqnarray}
\partial_t \mathbf{M}_1+\nabla\cdot \mathbf{M}_2 =\mathbf{F}+O(\Delta t),
\end{eqnarray}
\begin{eqnarray}
\partial_t \mathbf{M}_1+\nabla\cdot \mathbf{M}_2 + \frac{\Delta t}{\bar{A}_{10}}\big[\bar{A}_{20}\partial_t \mathbf{F}+\bar{B}_{10}\partial_t \mathbf{F}\big]\nonumber\\
= \mathbf{F} -\frac{\Delta t}{\bar{A}_{10}}\big[(\bar{A}_{21}-\bar{A}_{20})\partial_t \nabla\cdot \mathbf{M}_2 + \bar{B}_{11}\nabla\cdot \mathbf{M}_{2F}+\bar{A}_{22}\nabla\nabla : \mathbf{M}_3\big]
+O(\Delta t^2),
\end{eqnarray}
\end{subequations}
where
\begin{equation}\label{eq:ME-NSEs-3-1}
\partial_t^2 \mathbf{M}_1=-\partial_t\nabla\cdot \mathbf{M}_2 + \partial_t \mathbf{F}+O(\Delta t)
\end{equation}
is used. If we set
\begin{equation}\label{eq:M-NSEs-3}
\bar{A}_{20}\partial_t \mathbf{F}+\bar{B}_{10}\partial_t \mathbf{F}=(\bar{A}_{20}+\bar{B}_{10})\partial_t \mathbf{F}=0,
\end{equation}
the moment equation at the leading order of $O(\Delta t^2)$ can be derived from Eq. (\ref{eq:ME-NSEs-3}b),
\begin{eqnarray}\label{eq:ME-NSEs-3-2}
\partial_t \mathbf{M}_1+\nabla\cdot \mathbf{M}_2=\mathbf{F}-\frac{\Delta t}{\bar{A}_{10}}\big[(\bar{A}_{21}-\bar{A}_{20})\partial_t \nabla\cdot \mathbf{M}_2+\bar{B}_{11}\nabla\cdot \mathbf{M}_{2F}+\bar{A}_{22}\nabla\nabla : \mathbf{M}_3\big]+O(\Delta t^2).
\end{eqnarray}

Using Eq. (\ref{eq:ME-NSEs-3}a), and setting $-B_{11}\mathbf{F}=(A_{21}-A_{20})\mathbf{F}$, Eq. (\ref{eq:ME-NCDE-3-2-1}) becomes
\begin{eqnarray}\label{eq:ME-NSEs-0}
\partial_t M_0+\nabla\cdot \mathbf{M}_1=S-\frac{\Delta t}{A_{10}}\big[(A_{22}-(A_{21}-A_{20}))\nabla\nabla : \mathbf{M}_2\big]+O(\Delta t^2).
\end{eqnarray}

Furthermore, according to the NSEs, and choosing the moments and parameters in Eqs. (\ref{eq:ME-NCDE-3-2-1}) and (\ref{eq:ME-NSEs-3-2}) properly, we can correctly recover the NSEs from Eqs. (\ref{eq:ME-NCDE-3-2-1}) (or (\ref{eq:ME-NSEs-0})) and (\ref{eq:ME-NSEs-3-2}). For example, taking
\begin{subequations}\label{eq:M-NSEs-4}
\begin{equation}
M_0=\rho, \mathbf{M}_1=\rho\mathbf{ u},\mathbf{M}_2=\rho \mathbf{uu}+p\mathbf{I}, \mathbf{M}_3= c_s^2 \mathbf{\Delta} \cdot\rho \mathbf{u}+\mathbf{M}_{30},
\end{equation}
\begin{equation}
\mathbf{M}_{2F}=-\frac{1}{\bar{B}_{11}}\big[(\bar{A}_{21}-\bar{A}_{20})\partial_t \mathbf{M}_2+\bar{A}_{22}(c_s^2(\nabla\cdot\rho \mathbf{u I}+\mathbf{u}\nabla\rho+(\mathbf{u}\nabla\rho)^T)+ \nabla\cdot\mathbf{M}_{30})\big],
\end{equation}
\begin{equation}
\nu=-\frac{\bar{A}_{22}}{\bar{A}_{10}}\Delta t c_s^2
\end{equation}
\end{subequations}
in Eqs. (\ref{eq:ME-NSEs-0}) and (\ref{eq:ME-NSEs-3-2}) leads to the general NSEs,
\begin{subequations}\label{eq:NSEs}
\begin{equation}
\partial_t \rho+\nabla\cdot \rho \mathbf{u}=S-\frac{\Delta t}{A_{10}}\nabla\cdot\big[(A_{22}-(A_{21}-A_{20}))\nabla\cdot (\rho\mathbf{uu}+p\mathbf{I})\big]+O(\Delta t^2),
\end{equation}
\begin{equation}
\partial_t \rho \mathbf{u}+\nabla\cdot (\rho\mathbf{uu}+p\mathbf{I})=\mathbf{F}+\nu\nabla\cdot\rho\big[\nabla \mathbf{u}+(\nabla \mathbf{u})^T\big]+O(\Delta t^2),
\end{equation}
\end{subequations}
where the auxiliary source term $\mathbf{M}_{30}$ contains the additional term caused by the anisotropy of lattice tensor in a given lattice model, as mentioned in Ref. \cite{ChaiYuanShi2023}. If an isotropic lattice is used, one can take $\mathbf{M}_{30}=\mathbf{0}$.

For the incompressible flow, $\nabla\cdot (\rho\mathbf{uu}+p\mathbf{I})= O(Ma^2)$, thus Eq. (\ref{eq:NSEs}(a)) becomes
\begin{equation}\label{eq:NSEs-0}
\partial_t \rho+\nabla\cdot \rho \mathbf{u}=S+O(\Delta t^2+\Delta t Ma^2),
\end{equation}
therefore, the incompressible NSEs are recovered from Eqs. (\ref{eq:NSEs-0}) and (\ref{eq:NSEs}(b)) in the weak compressible limit, and $\mathbf{M}_{2F}$ can be simplified \cite{ChaiYuanShi2023}.

\subsection{The equilibrium, and source distribution functions of MesoDF method for NCDE and NSEs}

From above analysis, one can clearly observe that to recover the macroscopic PDEs, such as NCDE and NSEs from MesoDF method (\ref{eq:2-1}), the equilibrium, and source distribution functions should satisfy the constraints on basic moments, as depicted above. The common quadratic form of the EDF can be found in the previous work  \cite{ChaiYuanShi2023,LuChaiShi2011}, in which a general quadratic form on a rD$d$Q$q$ lattice was given by extending the commonly used one as 
\begin{equation}\label{eq:3-3}
g_j^{eq}(A,\mathbf{B},\mathbf{C})=\omega_j\left(A+\tilde{\textbf{c}}_j\cdot{\mathbf{B}}
                   +\tilde{\textbf{Q}}_j:\mathbf{C}\right),
\end{equation}
where $\omega_j$ is the weight coefficients (see Appendix A), and
\begin{equation}\label{eq:3-4}
\mathbf{\tilde{c}}_{j\alpha}=\mathbf{c}_{j\alpha}/c_{s\alpha}^2,
\mathbf{\tilde{Q}}_{j\alpha\alpha}=\mathbf{Q}_{j\alpha\alpha}/(c_{s\alpha}^2(c_\alpha^2-c_{s\alpha}^2)),
\mathbf{\tilde{Q}}_{j\alpha\beta}=\mathbf{Q}_{j\alpha\beta}/(2 c_{s\alpha}^2 c_{s\beta}^2) \ \ (\alpha \neq\beta),
\mathbf Q_j=\mathbf c_j \mathbf c_j -\mathbf \Delta^{(2)}.
\end{equation}
Here  $\mathbf{\Delta}^{(2)}=\sum_j \omega_j\mathbf{c}_j\mathbf{c}_j=\mathbf{diag}\{c_{s\alpha}\}$. We would also like to point out that strictly speaking, $c_{s\alpha}$ is direction-dependent, while it is usually taken to be a direction-independent form with $c_{s\alpha}=c_s$, $c_s$ is an adjustable parameter \cite{dHumieres1992,Ginzburg2005a,ChaiYuanShi2023}.

According to the properties of rD$d$Q$q$ lattice models, we can obtain the basic moments of $g_j^{eq}$ on the rD$d$Q$q$ lattice \cite{ChaiYuanShi2023,LuChaiShi2011},
\begin{subequations}\label{eq:3-5}
\begin{equation}
\sum_j g_j^{eq}=A,\ \ \sum_j \mathbf{c}_j g_j^{eq}=\mathbf{B},\ \ \sum_j \mathbf{c}_j \mathbf{c}_j
g_j^{eq}=A \mathbf{\Delta}^{(2)}+\mathbf{C},
\end{equation}
\begin{equation}
\sum_j \mathbf{c}_j \mathbf{c}_j\mathbf{c}_j g_j^{eq}=\mathbf{\Delta}^{(4)}\cdot \tilde{\mathbf{B}},
\end{equation}
\end{subequations}
where $\tilde{\mathbf{B}}_{\alpha}=\mathbf B_{\alpha}/ c_{s\alpha}^2$, and $\mathbf{\Delta}^{(4)}=\sum_j \omega_j\mathbf{c}_j\mathbf{c}_j\mathbf{c}_j\mathbf{c}_j=<\mathbf{\Delta}^{(2)}\mathbf{\Delta}^{(2)}>+\delta^{(4)}$ which is defined as
\begin{subequations}\label{eq:3-6}
\begin{equation}
<\mathbf{\Delta}^{(2)}\mathbf{\Delta}^{(2)}>_{\alpha\beta\gamma\theta}=\Delta_{\alpha\beta}^{(2)}\Delta_{\gamma\theta}^{(2)}+\Delta_{\alpha\gamma}^{(2)}\Delta_{\beta\theta}^{(2)}+\Delta_{\beta\gamma}^{(2)}\Delta_{\alpha\theta}^{(2)},\\
\end{equation}
\begin{equation}
\mathbf{\delta}^{(4)}_{\alpha\beta\gamma\theta}=c_{s\alpha}^2(c_\alpha^2-3c_{s\alpha}^2),\alpha=\beta=\gamma=\theta;\\
\mathbf{\delta}^{(4)}_{\alpha\beta\gamma\theta}=0, else.
\end{equation}
\end{subequations}

With the help of Eq. (\ref{eq:3-6}), we can rewrite Eq. (\ref{eq:3-5}b) as
\begin{equation}\label{eq:3-7}
\sum_j \mathbf{c}_{j\alpha} \mathbf{c}_{j\beta}\mathbf{c}_{j\gamma} g_j^{eq}=\Delta_{\alpha\beta}^{(2)}B_\gamma+\Delta_{\alpha\gamma}^{(2)}B_\beta+\Delta_{\beta\gamma}^{(2)}B_\alpha+\delta^{(4)}_{\alpha\beta\gamma\theta} \tilde{B}_\theta.
\end{equation}

For the NCDE, to satisfy the moment conditions in Eqs. (\ref{eq:M-NCDE-4}) and (\ref{eq:M-NCDE-6}), the explicit expressions of $f_{j}^{eq}$ and $F_j$ with $c_{s \alpha}=c_s$ for all $\alpha$ can be obtained by
\begin{eqnarray}\label{eq:NCDEfeq}
f_j^{eq} & = & g_j^{eq}(\phi, \mathbf{B}, \beta c_s^2\mathbf{D}+\mathbf{C}- c_s^2 \phi\mathbf{I}) \nonumber\\
         & = & \omega_j\left[\phi+\frac{\textbf{c}_{j\alpha}{\mathbf{B}_{\alpha}}}{c_s^2}
                   +\frac{(\beta c_s^2\mathbf{D}+\mathbf{C}- c_s^2 \phi\mathbf{I})_{\alpha\alpha}(c_{j\alpha}^2 - c_s^2)}{c_s^2(c_{\alpha}^2-c_s^2)}+\frac{(\beta c_s^2\mathbf{D}+\mathbf{C}-  c_s^2 \phi\mathbf{I})_{\alpha\bar{\alpha}}(c_{j\alpha}c_{j\bar{\alpha}})}{2c_s^4}\right],
\end{eqnarray}
\begin{equation}\label{eq:SourceF}
F_j=g_j^{eq}(S,\mathbf{M}_{1F},\mathbf{0})=\omega_j (S+\frac{\mathbf{c}_{j\alpha} \mathbf{M}_{1F,\alpha}}{c_s^2}),
\end{equation}
where $\bar{\alpha}$ denotes the index $\gamma$ with $\gamma\neq\alpha$, $\mathbf{M}_{1F}$ is given by Eq. (\ref{eq:M-NCDE-6}).

While for the NSEs, to satisfy the moment conditions in Eqs. (\ref{eq:M-NSEs-2}) and (\ref{eq:M-NSEs-4}), $f_{j}^{eq}$ and $F_j$ with $c_{s \alpha}=c_s$ for all $\alpha$ can be given by
\begin{eqnarray}\label{eq:NSEfeq}
f_j^{eq} & = & g_j^{eq}(\rho, \rho\mathbf{u}, \rho\mathbf{uu}) = \rho g_j^{eq}(1, \mathbf{u}, \mathbf{uu})\nonumber\\
         & = & \omega_j \rho \left[1+\frac{\textbf{c}_{j\alpha}{\mathbf{u}_{\alpha}}}{c_s^2}
                   +\frac{u_{\alpha}^2(c_{j\alpha}^2 - c_s^2)}{c_s^2(c_{\alpha}^2-c_s^2)}+\frac{u_{\alpha}u_{\bar{\alpha}}(c_{j\alpha}c_{j\bar{\alpha}})}{2c_s^4}\right],
\end{eqnarray}
\begin{equation}\label{eq:NSEFeq}
F_j=g_j^{eq}(S, \mathbf{F}, \mathbf{M}_{2F}-c_s^2 S\mathbf{I})=\omega_j \left[S+\frac{\textbf{c}_{j\alpha} {\mathbf{\mathbf{F}}_{\alpha}}}{c_s^2}
                   +\frac{(\mathbf{M}_{2F,\alpha\alpha}-c_s^2 S)(c_{j\alpha}^2 - c_s^2)}{c_s^2(c_{\alpha}^2-c_s^2)}+\frac{\mathbf{M}_{2F,\alpha\bar{\alpha}}(c_{j\alpha}c_{j\bar{\alpha}})}{2c_s^4}\right],
\end{equation}
where $\mathbf{M}_{2F}$ is determined by Eq. (\ref{eq:M-NSEs-4}b).

\textbf{Remark 4} : For  incompressible flow, inserting $\rho=\rho_0+\delta \rho$ with $\delta \rho=O(Ma^2)$ into Eq. (\ref{eq:NSEfeq}), and omitting the terms of $O(Ma^3)$, one can obtain the He-Luo model \cite{HeLuo1997} of the present MesoFD scheme with $f_j^{eq}=g_j^{eq}(\rho, \rho_0\mathbf{u}, \rho_0\mathbf{uu})$.

\section{Numerical tests}
In this section, several different numerical examples are used to test the present MesoFD scheme, while for simplicity, we only consider the MMFD scheme (\ref{eq:3layerMMFD-1}).
To test the accuracy of the MesoFD scheme, the following global relative errors (GRE) are adopted,

$$
\operatorname{GRE}=\frac{\sum_i\left|\phi\left(\boldsymbol{x}_i, t\right)-\phi^*\left(\boldsymbol{x}_i, t\right)\right|}{\sum_i\left|\phi^*\left(\boldsymbol{x}_i, t\right)\right|},
$$
where $\phi(\boldsymbol{x}, t)$ and $\phi^*(\boldsymbol{x}, t)$ denote numerical and analytical solutions, respectively. 

{\bf{Example 1}}: Isotropic convection-diffusion equation

The following two-dimensional isotropic CDE is first considered to test the accuracy of the MesoFD scheme,

$$
\partial_t \phi+u_x \partial_x \phi+u_y \partial_y \phi=\kappa\left(\partial_{x x} \phi+\partial_{y y} \phi\right)+R,
$$
where $u_x=u_y=0.1, ~\kappa$ is the diffusion coefficient, and the source term $R$ is given by

$$
\begin{aligned}
R=  \exp \left[\left(1-2 \pi^2 \kappa\right) t\right]\left( \sin [\pi(x+y)]  +\pi\left(u_x+u_y\right) \cos [\pi(x+y)]\right).
\end{aligned}
$$

The analytical solution of this problem can be expressed as

$$
\phi(x, y, t)=\exp \left[\left(1-2 \pi^2 \kappa\right) t\right] \sin (\pi(x+y)).
$$
The computational domain is $(x,y)\in[0,2]\times[0,2]$ and $t\in[0,T]$.

We take 
\begin{equation*}
f_j^{eq}=\omega_j\phi \left[1+\frac{\mathbf{c}_j\cdot\mathbf{u}}{c_s^2}+\frac{(\mathbf{u}\mathbf{u}+\beta c_s^2\mathbf{I}-c_s^2\mathbf{I}):(\mathbf{c}_j\mathbf{c}_j-c_s^2\mathbf{I})}{2c_s^4}\right],
\end{equation*} 
and
\begin{equation*}
F_j =\omega_j R \left[1+\frac{\mathbf{c}_j\cdot\mathbf{u}}{c_s^2}\right].
\end{equation*} 
In the simulations, we adopt the D2Q9 lattice model with $\omega_0 = 4/9,~\omega_{1-4} = 1/9,~\omega_{5-8} = 1/36$, and $c_s^2=c^2/3,~c=\Delta x/ \Delta t$. We consider the following cases:

\textbf{Case 1}: SRT LB model with $\tau =1$. Let $A_{10}=A_{11}=-B_{00}=-1,~A_{21}=0,~A_{22}=\frac{1}{2},~a_{0}=0$, then we have \begin{align*}
    a_{-1}=1,\quad a_0=0,\quad a_1=0,\quad \hat{a}_{-1}=0,\quad \hat{a}_0=0,\quad \hat{a}_1=0, \\
    b_0=b_1=\hat{b}_0=\hat{b}_1=0,\quad b_{-1}=\frac{3}{2},\quad \hat{b}_{-1}=-\frac{1}{2}, 
  \end{align*}
  the MMFD scheme (\ref{eq:3layerMMFD-1}) becomes
  \begin{align*}
    \phi (\mathbf{x},t+\Delta t)=\sum_j\left(f_j^{eq}(\mathbf{x}-\mathbf{c}_j\Delta t,t)+\frac{3}{2}\Delta t F_j (\mathbf{x}-\mathbf{c}_j\Delta t,t)-\frac{1}{2}\Delta t F_j (\mathbf{x}-\mathbf{c}_j\Delta t,t-\Delta t)\right).
  \end{align*}
  
  \textbf{Case 2}: Let $A_{10}=A_{11}=-B_{00}=-1,~A_{21}=0,~A_{22}=\frac{1}{2},~a_{0}=\frac{1}{4}$, and we can derive
  \begin{align*}
  a_{-1}=\frac{7}{8},\quad a_0=\frac{1}{4},\quad a_1=-\frac{1}{8},\quad \hat{a}_{-1}=\frac{1}{8},\quad \hat{a}_0=-\frac{1}{4},\quad \hat{a}_1=-\frac{1}{8}, \\
  b_0=b_1=\hat{b}_0=\hat{b}_1=0,\quad b_{-1}=\frac{3}{2},\quad \hat{b}_{-1}=-\frac{1}{2}, 
  \end{align*}
 and 
 \begin{align*}
	  \phi(\mathbf{x},t+\Delta t)&=\sum_j\left( \frac{7}{8}f_j^{eq}(\mathbf{x}-\mathbf{c}_j\Delta t,t)+\frac{1}{4}f_j^{eq}(\mathbf{x},t)-\frac{1}{8}f_j^{eq}(\mathbf{x}+\mathbf{c}_j\Delta t,t) \right. \\
	  &+\frac{1}{8}f_j^{eq}(\mathbf{x}-\mathbf{c}_j\Delta t,t-\Delta t)-\frac{1}{4}f_j^{eq}(\mathbf{x},t-\Delta t)-\frac{1}{8}f_j^{eq}(\mathbf{x}+\mathbf{c}_j\Delta t,t-\Delta t) \\
	  &\left. +\frac{3}{2}\Delta t F_j (\mathbf{x}-\mathbf{c}_j\Delta t,t)-\frac{1}{2}\Delta t F_j (\mathbf{x}-\mathbf{c}_j\Delta t,t-\Delta t)\right).  
 \end{align*}
  
  \textbf{Case 3}: Let $A_{10}=A_{11}=-B_{00}=-1,~A_{21}=\frac{1}{4},~A_{22}=\frac{1}{2},~a_{0}=0$, then we have
  \begin{align*}
  a_{-1}=\frac{7}{8},\quad a_0=0,\quad a_1=\frac{1}{8},\quad \hat{a}_{-1}=\frac{1}{8},\quad \hat{a}_0=0,\quad \hat{a}_1=-\frac{1}{8}, \\
  b_0=b_1=\hat{b}_0=\hat{b}_1=0,\quad b_{-1}=\frac{3}{2},\quad \hat{b}_{-1}=-\frac{1}{2}, 
  \end{align*}
and
  \begin{align*}
  \phi(\mathbf{x},t+\Delta t)&=\sum_j\left( \frac{7}{8}f_j^{eq}(\mathbf{x}-\mathbf{c}_j\Delta t,t)+\frac{1}{8}f_j^{eq}(\mathbf{x}+\mathbf{c}_j\Delta t,t)+\frac{1}{8}f_j^{eq}(\mathbf{x}-\mathbf{c}_j\Delta t,t-\Delta t)-\frac{1}{8}f_j^{eq}(\mathbf{x}+\mathbf{c}_j\Delta t,t-\Delta t) \right. \\
  &\left.+\frac{3}{2}\Delta t F_j (\mathbf{x}-\mathbf{c}_j\Delta t,t)-\frac{1}{2}\Delta t F_j (\mathbf{x}-\mathbf{c}_j\Delta t,t-\Delta t)\right).  
  \end{align*}
  
  Let $N_x = \frac{1}{\Delta x}$, $N_t = \frac{1}{\Delta t}$, then the errors and convergency orders at $T=1.0$ are listed in Table \ref{table1}. From the table, the second-order convergency can be found, which is in agreement with the theoretical analysis.

  \begin{table}[H]
	\centering
	\caption{The errors and convergency orders for Case 1 - 3}
	\label{table1}
	\begin{tabular}{cccccccc}
		\toprule[1.5pt]
		\multicolumn{2}{l}{} & \multicolumn{2}{c}{Case 1} & \multicolumn{2}{c}{Case 2} &  \multicolumn{2}{c}{Case 3}   \\ \cline{3-8} 
		$N_x$        & $N_t$         & GRE1              & order1          & GRE2              & order2          & GRE3              & order3                     \\ \bottomrule[1.5pt]
		10        & 6          & 1.51E-01          &                 & 1.05E-01          &                 & 1.01E-01          &                 \\
  			20        & 24         & 4.16E-02          & 1.86            & 2.50E-02          & 2.07            & 2.45E-02          & 2.04            \\
  			40        & 96         & 1.07E-02          & 1.97            & 6.18E-03          & 2.01            & 6.08E-03          & 2.01            \\
  			80        & 384        & 2.68E-03          & 1.99            & 1.54E-03          & 2.00            & 1.52E-03          & 2.00            \\ \bottomrule[1.5pt]
	\end{tabular}
	\vspace{1em}
\end{table}

{\bf{Example 2}}: Isotropic convection-diffusion equation without the source term

For the homogeneous isotropic convection-diffusion equation,
\begin{equation*}
\partial_t \phi+u_x \partial_x \phi+u_y \partial_y \phi=D\left(\partial_{x x} \phi+\partial_{y y} \phi\right),
\end{equation*}
where $u_x=u_y=1.0$ and $D=0.02$ is the diffusion coefficient. The computational domain is $(x,y)\in [0,2]\times[0,2]$ and $t\in [0,T]$, and the exact solution can be expressed as 
$$\phi(x, y, t)=\exp \left(-2D \pi^2  t\right) \cos [\pi(x+y)-\pi (u_x+u_y)t].$$
Similar to above discussion, we also take
\begin{equation*}
f_j^{eq}=\omega_j\phi [1+\frac{\mathbf{c}_j\cdot\mathbf{u}}{c_s^2}+\frac{(\mathbf{u}\mathbf{u}+\beta c_s^2\mathbf{I}-c_s^2\mathbf{I}):(\mathbf{c}_j\mathbf{c}_j-c_s^2\mathbf{I})}{2c_s^4}],
\end{equation*}
and use D2Q9 lattice model with  $\omega_0 = 4/9,~\omega_{1-4} = 1/9,~\omega_{5-8} = 1/36$. 
Here, we also consider the following three cases:

\textbf{Case 1}: Let $A_{10}=A_{11}=-1,~A_{21}=0,~A_{22}=\frac{1}{2},~a_{0}=0,~\beta=1$, then we have
  \begin{align*}
    a_{-1}=1,\quad a_0=0,\quad a_1=0,\quad \hat{a}_{-1}=0,\quad \hat{a}_0=0,\quad \hat{a}_1=0,
  \end{align*}
  and
  \begin{align*}
    \phi (\mathbf{x},t+\Delta t)=\sum_j f_j^{eq}(\mathbf{x}-\mathbf{c}_j\Delta t,t).
  \end{align*}
  
  \textbf{Case 2}: Let $A_{10}=A_{11}=-1, ~A_{21}=\frac{1}{4}, ~A_{22}=\frac{1}{2},~a_{0}=0,~\beta=1$, then we obtain
  \begin{align*}
  a_{-1}=\frac{7}{8},\quad a_0=0,\quad a_1=\frac{1}{8}, \quad \hat{a}_{-1}=\frac{1}{8},\quad \hat{a}_0=0,\quad \hat{a}_1=-\frac{1}{8},
  \end{align*}
  and
  \begin{align*}
	  \phi(\mathbf{x},t+\Delta t)=\sum_j\left(\frac{7}{8}f_j^{eq}(\mathbf{x}-\mathbf{c}_j\Delta t,t)+\frac{1}{8}f_j^{eq}(\mathbf{x}+\mathbf{c}_j\Delta t,t)
	  +\frac{1}{8}f_j^{eq}(\mathbf{x}-\mathbf{c}_j\Delta t,t-\Delta t)-\frac{1}{8}f_j^{eq}(\mathbf{x}+\mathbf{c}_j\Delta t,t-\Delta t)\right).  
  \end{align*}
  
  \textbf{Case 3}: Let $A_{10}=A_{11}=-1,~A_{21}=0,~A_{22}=1,~a_{0}=0,~\beta=\frac{1}{2}$, then we derive
  \begin{align*}
  a_{-1}=1,\quad a_0=a_1=0,\quad \hat{a}_{-1}=\frac{1}{2},\quad \hat{a}_0=-1\quad \hat{a}_1=\frac{1}{2},
  \end{align*}
  and
  \begin{align*}
  \phi (\mathbf{x},t+\Delta t)=\sum_j\left(f_j^{eq}(\mathbf{x}-\mathbf{c}_j\Delta t,t)+\frac{1}{2}f_j^{eq}(\mathbf{x}+\mathbf{c}_j\Delta t,t)-f_j^{eq}(\mathbf{x}-\mathbf{c}_j\Delta t,t-\Delta t)+\frac{1}{2}f_j^{eq}(\mathbf{x}+\mathbf{c}_j\Delta t,t-\Delta t)\right).  
  \end{align*}

  \begin{table}[H]
	\centering
	\caption{The errors and convergency orders for Case 1 - 3}
	\label{table2}
	\begin{tabular}{cccccccc}
		\toprule[1.5pt]
		\multicolumn{2}{l}{} & \multicolumn{2}{c}{Case 1} & \multicolumn{2}{c}{Case 2} &  \multicolumn{2}{c}{Case 3}   \\ \cline{3-8} 
		$N_x$        & $N_t$         & GRE1              & order1          & GRE2              & order2          & GRE3              & order3                     \\ \bottomrule[1.5pt]
		10        & 12          & 5.44E-01          &                 & 5.74E-01          &                 & 7.39E-00         &                 \\
  			20        & 48         & 7.20E-02          & 2.92            & 1.41E-02          & 2.02            & 2.41E-02          & 4.94            \\
  			40        & 192         & 1.45E-02          & 2.31            & 3.31E-02          & 2.10            & 6.24E-02          & 1.95            \\
  			80        & 768        & 3.34E-03          & 2.12            & 8.05E-03          & 2.04            & 1.58E-02          & 1.98            \\ \bottomrule[1.5pt]
	\end{tabular}
	\vspace{1em}
\end{table}

 In our simulations, $N_x = \frac{1}{\Delta x}$, $N_t = \frac{1}{\Delta t}$, and the errors and convergency orders at $T=1.0$ are measured in Table \ref{table2}. As shown in this table, all of the three cases are of a second-order convergency also can be found, which agrees well with the theoretical analysis. In order to test the stable condition of the two-level scheme, i.e.,
\begin{align*}
    \frac{\mathbf{u}^2}{c_s^2}\leq \frac{1}{1-\omega_0}=\frac{9}{5},
\end{align*}
the case 1 is adopted with $D=0.01,~u_x=u_y=1.0$.  We calculate the errors, and  list them in Table \ref{table3}. From the tabel, we can see that if $\frac{\mathbf{u}^2}{c_s^2}$ is more than 9/5, the scheme is unstable which is in agreement with Theorem 1.

  \begin{table}[H]
	\centering
\caption{The errors under different values of $c_s^2$}
	\label{table3}
{\begin{tabular}{cccc}
\toprule
$\Delta x$  & $\Delta t$   & $\frac{\mathbf{u}^2}{c_s^2}$  & GRE   \\ 
\midrule
$\frac{1}{10}$ & $\frac{1}{6}$  & $\frac{50}{3}$ & 7.78E+02 \\
$\frac{1}{20}$ & $\frac{1}{24}$ & $\frac{25}{6}$ & 2.02E-00  \\
$\frac{1}{40}$ & $\frac{1}{96}$ & $\frac{25}{24}$ & 3.26E-02 \\
$\frac{1}{80}$ & $\frac{1}{384}$ & $\frac{25}{96}$ & 6.90E-03 \\ 
\bottomrule
\end{tabular}
}
\end{table}

{\bf Example 3:} Wave equation

 For the two-dimensional wave equation
  \begin{equation*}
    \frac{\partial^2 \phi}{\partial t^2}=\frac{1}{2}(\frac{\partial^2 \phi}{\partial x^2}+\frac{\partial^2 \phi}{\partial y^2}), \quad (x,y)\in [0,1]\times [0,1],t\in [0,1],
  \end{equation*}
  the analytical solution can be given by
  $$\phi (x,y,t)=e^{x+y+t}.$$
  In the following simulations, we take $N_x=N_y$, and consider D2Q5 and D2Q9 lattice models. Based on the idea of the finite difference method, we have the following approximation to the wave equation,
  \begin{align*}
    \frac{2}{\Delta t}\left[\frac{\phi (\mathbf{x},\Delta t)-\phi (\mathbf{x},0)}{\Delta t}-\frac{\partial \phi (\mathbf{x},0)}{\partial t}\right]=\frac{1}{2}\sum_{j=1}^2\left[\frac{\phi (\mathbf{x}-\mathbf{c}_j\Delta t,0)-2\phi (\mathbf{x},0)+\phi (\mathbf{x}+\mathbf{c}_j\Delta t,0)}{\Delta x^2}\right],
  \end{align*}
  that is,
  \begin{align*}
    \phi (\mathbf{x},\Delta t)=\phi (\mathbf{x},0)+\Delta t\frac{\partial \phi (\mathbf{x},0)}{\partial t}+\left(\frac{\Delta t}{2\Delta x}\right)^2 \sum_{j=1}^2[\phi (\mathbf{x}-\mathbf{c}_j\Delta t,0)-2\phi (\mathbf{x},0)+\phi (\mathbf{x}+\mathbf{c}_j\Delta t,0)].
  \end{align*}

  The equilibrium function is taken as follows
  \begin{equation*}
    f_j^{eq}=\omega_j\phi \left[1+\frac{(\beta c_s^2\mathbf{I}-c_s^2\mathbf{I}):(\mathbf{c}_j\mathbf{c}_j-c_s^2\mathbf{I})}{2c_s^4}\right]=\omega_j\phi \left[1+\frac{(\beta -1)(|\mathbf{c}_j|^2-2c_s^2)}{2c_s^2}\right].
  \end{equation*}
  Here $a_0=2(1-\gamma),~a_1=a_{-1}=\gamma,~\hat{a}_0=-1,~\hat{a}_1=\hat{a}_{-1}=0$, $A_{22}=\gamma$. 

  \textbf{Case 1}: Take the D2Q5 model and let $\beta=1,\gamma=\frac{3}{8\beta}=\frac{3}{8}$, we have
  \begin{align*}
  \phi (\mathbf{x},t+\Delta t)&=\sum_j\left(\frac{3}{8}f_j^{eq}(\mathbf{x}-\mathbf{c}_j\Delta t,t)+\frac{5}{4}f_j^{eq}(\mathbf{x},t)+\frac{3}{8}f_j^{eq}(\mathbf{x}+\mathbf{c}_j\Delta t,t)-f_j^{eq}(\mathbf{x},t-\Delta t)\right).
  \end{align*}

  \textbf{Case 2}: Consider the D2Q9 model and the parameters are the same as in Case 1.

  \textbf{Case 3}: Consider the D2Q9 model and take $\beta=\frac{3}{4},~\gamma=\frac{3}{8\beta}=\frac{1}{2}$, we obtain
  \begin{align*}
  \phi (\mathbf{x},t+\Delta t)&=\sum_j\left(\frac{1}{2}f_j^{eq}(\mathbf{x}-\mathbf{c}_j\Delta t,t)+f_j^{eq}(\mathbf{x},t)+\frac{1}{2}f_j^{eq}(\mathbf{x}+\mathbf{c}_j\Delta t,t)-f_j^{eq}(\mathbf{x},t-\Delta t)\right).
  \end{align*}

It is noted that the scheme in Case 1 is  an explicit finite-difference scheme. We measure the errors and convergency orders  in Table \ref{table4}, and from which second-order convergency can be found. 

  \begin{table}[H]
	\centering
	\caption{The errors and convergency orders of Case 1-3 for the wave equation}
	\label{table4}
	\begin{tabular}{cccccccc}
		\toprule[1.5pt]
		\multicolumn{2}{l}{} & \multicolumn{2}{c}{Case 1} & \multicolumn{2}{c}{Case 2} &  \multicolumn{2}{c}{Case 3}   \\ \cline{3-8} 
		$N_x$  & $N_t$ & GRE1 & order & GRE2 & order2 & GRE3 & order3                     \\ \bottomrule[1.5pt]
     5 & 10 & 2.86E-04 &      & 6.04E-03 &      & 3.91E-03 &   \\
    10 & 20 & 6.09E-05 & 2.23 & 1.74E-03 & 1.79 & 1.13E-03 & 1.79  \\
    20 & 40 & 1.39E-05 & 2.13 & 4.55E-04 & 1.94 & 2.97E-04 & 1.93  \\
    40 & 80 & 3.30E-06 & 2.07 & 1.16E-04 & 1.98 & 7.58E-05 & 1.97  \\ \bottomrule[1.5pt]
	\end{tabular}
	\vspace{1em}
\end{table}


\section{Conclusions}

By utilizing the equilibrium distribution function to represent the distribution function and/or the nonequilibrium distribution function, a new kind of the mesoscopic numerical method named the mesoscopic finite difference method, with a multi-level and multi-point structure is presented for some PDEs, including but not limited to the NSEs and NCDE. The macroscopic version of the MesoFD scheme, the MMFD scheme for the conserved variables, can be obtained directly from the moments of the MesoFD scheme.

A detailed DTE analysis on the present MesoFDM is conducted at the acoustic scaling, and the macroscopic moment equations (MEs) with truncation errors of different orders can be obtained. Based on the MEs with second- and third-order truncation errors, the NSEs, NCDE and wave equation can be recovered with the second-order accuracy in space. If the MEs with a higher-order truncation errors are used, one can recover the related PDEs with the higher-order accuracy, and more lattice velocities are needed. The MMFD scheme is a unified FD scheme for PDEs in any dimensions. The existing FD schemes derived from the LB models and even the traditional FD schemes are the special cases of present MMFD scheme. The MesoFDM  not only establishes a connection between the MesoNM and the traditional macroscopic FDM, but also extends the LB method and traditional FDM. We can combine the 'bottom-up' mesoscopic method with the 'top-down' macroscopic method to design more effective numerical methods by using their own advantages.

\section*{Acknowledgements}
This work was financially supported by the National Natural Science Foundation of China (Grants No. 12072127 and No. 12202130) and Jiangsu Provincial Scientific
Research Center of Applied Mathematics, China (No. BK20233002).

\appendix
\section{\label{app:sec1}Lattice models of RMRT-LB method on rectangular lattice \cite{ChaiYuanShi2023}}

In this appendix, we give some commonly used rD$d$Q$q$ lattice models in Ref. \cite{ChaiYuanShi2023}. Let $c_\alpha=\Delta x_{\alpha}/\Delta t, d_{0\alpha}=c^2_{s\alpha}/c^2_{\alpha}$ ($\alpha=1, 2, \ldots, d$) in $d$-dimensional space with $\Delta x_{\alpha}$ being the spacing step in $\alpha$ axis. In this case, the discrete velocities and weight coefficients in the common rD$d$Q$q$ lattice models can be given as follows.

rD2Q9 lattice:

\begin{equation}\label{eq:D2Q9}
\begin{split}
&\{\mathbf{c}_j,0\leq j \leq 8\}=
\left(
\begin{array}{ccccccccc}
    0 & c_1 &   0 & -c_1  &    0 & c_1 & -c_1 & -c_1 &  c_1 \\
    0 &   0 & c_2 &    0  & -c_2 & c_2 &  c_2 & -c_2 & -c_2 \\
\end{array}
\right),\\
&\omega_j\geq 0, \omega_1=\omega_3,\omega_2=\omega_4,\omega_5=\omega_6=\omega_7=\omega_8, \omega_0=1-2(\omega_1+\omega_2+2\omega_5),\\
&\omega_0=(1-d_{01})(1-d_{02}), \omega_1=d_{01}(1-d_{02})/2, \omega_2=(1-d_{01})d_{02}/2, \omega_5=d_{01}d_{02}/4,
\end{split}
\end{equation}
where $d_{01},d_{02} \in (0,1)$.

rD2Q5I lattice:

\begin{equation}\label{eq:D2Q5I}
\begin{split}
&\{\mathbf{c}_j,0\leq j \leq 4\}=
\left(
\begin{array}{ccccc}
    0 & c_1 &   0 & -c_1  &    0  \\
    0 &   0 & c_2 &    0  & -c_2  \\
\end{array}
\right),\\
&\omega_j\geq 0, \omega_1=\omega_3,\omega_2=\omega_4, \omega_0=1-2(\omega_1+\omega_2),\\
& \omega_0=1-d_{01}-d_{02}, \omega_1=d_{01}/2, \omega_2=d_{02}/2,
\end{split}
\end{equation}
where $d_{01}>0,d_{02}>0$, and $d_{01}+d_{02}<1$.

rD2Q5II lattice:

\begin{equation}\label{eq:D2Q5II}
\begin{split}
&\{\mathbf{c}_j,0\leq j \leq 4\}=
\left(
\begin{array}{ccccc}
    0 & c_1 & -c_1 & -c_1 &  c_1 \\
    0 & c_2 &  c_2 & -c_2 & -c_2 \\
\end{array}
\right),\\
&\omega_j\geq 0, \omega_1=\omega_2=\omega_3=\omega_4, \omega_0=1-4\omega_1,\\
& \omega_0=1-d_0, \omega_1=d_0/4,
\end{split}
\end{equation}
where $d_{01}=d_{02}=d_0 \in (0,1)$.

rD3Q27 lattice:

\begin{subequations}\label{eq:D3Q27}
\begin{equation}\label{eq:D3Q27a}
\{\mathbf{c}_j,0\leq j \leq 6\}=
\left(
\begin{array}{ccccccc}
    0 & c_1 &   0 &   0 & -c_1  &    0 &    0 \\
    0 &   0 & c_2 &   0 &    0  & -c_2 &    0 \\
    0 &   0 &   0 & c_3 &    0  &    0 & -c_3 \\
\end{array}
\right),
\end{equation}
\begin{equation}\label{eq:D3Q27b}
\{\mathbf{c}_j,7\leq j \leq 18\}=
\left(
\begin{array}{cccccccccccc}
    c_1 & -c_1 & -c_1 &  c_1 & c_1 & -c_1 & -c_1 & c_1 &   0 &    0 &    0 &    0 \\
    c_2 &  c_2 & -c_2 & -c_2 &   0 &    0 &    0 &   0 & c_2 & -c_2 & -c_2 &  c_2 \\
      0 &    0 &    0 &    0 & c_3 & -c_3 & -c_3 & c_3 & c_3 &  c_3 & -c_3 & -c_3 \\
\end{array}
\right),
\end{equation}
\begin{equation}\label{eq:D3Q27c}
\{\mathbf{c}_j,19\leq j \leq 26\}=
\left(
\begin{array}{cccccccc}
    c_1 &  c_1 &  c_1 & -c_1 & -c_1 & -c_1 & -c_1 &  c_1     \\
    c_2 &  c_2 & -c_2 &  c_2 & -c_2 & -c_2 &  c_2 & -c_2     \\
    c_3 & -c_3 &  c_3 &  c_3 & -c_3 &  c_3 & -c_3 & -c_3     \\
\end{array}
\right),
\end{equation}
\begin{equation*}\label{eq:D3Q27d}
\begin{split}
\omega_j\geq 0, \omega_1&=\omega_4,\omega_2=\omega_5,\omega_3=\omega_6,\omega_7=\omega_8=\omega_9=\omega_{10}, \omega_{11}=\omega_{12}=\omega_{13}=\omega_{14},\\
&\omega_{15}=\omega_{16}=\omega_{17}=\omega_{18},\omega_j=\omega_{19} (j>19), \omega_0=1-\sum_{j>0} \omega_j,
\end{split}
\end{equation*}
\end{subequations}
where
\begin{subequations}\label{eq:D3Q27W}
\begin{equation}
\omega_0=(1-d_{01})(1-d_{02})(1-d_{03}),
\end{equation}
\begin{equation}
\omega_1=d_{01}(1-d_{02})(1-d_{03})/2, \omega_2=d_{02}(1-d_{01})(1-d_{03})/2, \omega_3=d_{03}(1-d_{01})(1-d_{02})/2,
\end{equation}
\begin{equation}
\omega_7=d_{01}d_{02}(1-d_{03})/4, \omega_{11}=d_{01}d_{03}(1-d_{02})/4, \omega_{15}=d_{02}d_{03}(1-d_{01})/4,\omega_{19}=d_{01}d_{02}d_{03}/8,
\end{equation}
\end{subequations}
where $d_{01}, d_{02}, d_{03} \in (0,1) $.

rD3Q19 lattice:

\begin{subequations}\label{eq:D3Q19}
\begin{equation}\label{eq:D3Q19a}
\{\mathbf{c}_j,0\leq j \leq 6\}=
\left(
\begin{array}{ccccccc}
    0 & c_1 &   0 &   0 & -c_1  &    0 &    0 \\
    0 &   0 & c_2 &   0 &    0  & -c_2 &    0 \\
    0 &   0 &   0 & c_3 &    0  &    0 & -c_3 \\
\end{array}
\right),
\end{equation}
\begin{equation}\label{eq:D3Q19b}
\{\mathbf{c}_j,7\leq j \leq 18\}=
\left(
\begin{array}{cccccccccccc}
    c_1 & -c_1 & -c_1 &  c_1 & c_1 & -c_1 & -c_1 & c_1 &   0 &    0 &    0 &    0 \\
    c_2 &  c_2 & -c_2 & -c_2 &   0 &    0 &    0 &   0 & c_2 & -c_2 & -c_2 &  c_2 \\
      0 &    0 &    0 &    0 & c_3 & -c_3 & -c_3 & c_3 & c_3 &  c_3 & -c_3 & -c_3 \\
\end{array}
\right),
\end{equation}
\begin{equation*}\label{eq:D3Q27c}
\begin{split}
\omega_j\geq 0, \omega_1&=\omega_4,\omega_2=\omega_5,\omega_3=\omega_6,\omega_7=\omega_8=\omega_9=\omega_{10}, \omega_{11}=\omega_{12}=\omega_{13}=\omega_{14},\\
&\omega_{15}=\omega_{16}=\omega_{17}=\omega_{18}, \omega_0=1-\sum_{j>0} \omega_j,
\end{split}
\end{equation*}
\end{subequations}
where
\begin{subequations}\label{eq:D3Q19W}
\begin{equation}
\omega_0=1-d_{01}-d_{02}-d_{03}+d_{01}d_{02}+d_{01}d_{03}+d_{02}d_{03},
\end{equation}
\begin{equation}
\omega_1=d_{01}(1-d_{02}-d_{03})/2, \omega_2=d_{02}(1-d_{01}-d_{03})/2, \omega_3=d_{03}(1-d_{01}-d_{02})/2,
\end{equation}
\begin{equation}
\omega_7=d_{01}d_{02}/4, \omega_{11}=d_{01}d_{03}/4, \omega_{15}=d_{02}d_{03}/4,
\end{equation}
\end{subequations}
where $d_{01}, d_{02}$ and $d_{03} $ make the weights $\omega_j>0, \forall j$.

rD3Q9 lattice:

\begin{equation}\label{eq:D3Q9}
\begin{split}
&\{\mathbf{c}_j,0\leq j \leq 8\}=
\left(
\begin{array}{ccccccccc}
   0 & c_1 &  c_1 &  c_1 & -c_1 & -c_1 & -c_1 & -c_1 &  c_1     \\
   0 & c_2 &  c_2 & -c_2 &  c_2 & -c_2 & -c_2 &  c_2 & -c_2     \\
   0 & c_3 & -c_3 &  c_3 &  c_3 & -c_3 &  c_3 & -c_3 & -c_3     \\
\end{array}
\right),\\
&\omega_j\geq 0, \omega_1=\omega_2=...=\omega_8, \omega_0=1-8\omega_1,\\
& \omega_0=1-d_0, \omega_1=d_0/8,
\end{split}
\end{equation}
where $d_{01}=d_{02}=d_{03}=d_0 \in (0,1)  $ .

rD3Q7 lattice:

\begin{subequations}\label{eq:D3Q7}
\begin{equation}\label{eq:D3Q7a}
\{\mathbf{c}_j,0\leq j \leq 6\}=
\left(
\begin{array}{ccccccc}
    0 & c_1 &   0 &   0 & -c_1  &    0 &    0 \\
    0 &   0 & c_2 &   0 &    0  & -c_2 &    0 \\
    0 &   0 &   0 & c_3 &    0  &    0 & -c_3 \\
\end{array}
\right),
\end{equation}
\begin{equation*}\label{eq:D3Q7b}
\begin{split}
\omega_j\geq 0, \omega_1&=\omega_4,\omega_2=\omega_5,\omega_3=\omega_6, \omega_0=1-2(\omega_1+\omega_2+\omega_3),
\end{split}
\end{equation*}
\end{subequations}
where
\begin{equation}\label{eq:D3Q7W}
\omega_0=1-d_{01}-d_{02}-d_{03},\omega_1=d_{01}/2, \omega_2=d_{02}/2, \omega_3=d_{03}/2,
\end{equation}
where $d_{01}, d_{02}, d_{03} >0, d_{01} + d_{02} + d_{03}<1 $ .

\textbf{Remark} \uppercase\expandafter{\romannumeral1}: If we take $\omega_0=0$ or remove the rest velocity $\mathbf{c}_0=\mathbf{0}$ in the rD$d$Q$q$ lattice model mentioned above, the rD$d$Q$(q-1)$ lattice model can be obtained.

\bibliography{MesoLB}

\end{document}